\documentclass{amsart}
\usepackage{amscd,amsmath}
\usepackage[dvips]{graphicx}
\usepackage{hyperref}
\usepackage[latin1]{inputenc}
\usepackage{array}
\usepackage[all]{xy}

\newcommand{\sing}{{\rm{Sing}}}
\newcommand{\f}{{\mathcal F}}
\newcommand{\GF}{{G_{\f}}}
\newcommand{\rank}{{\rm{rank}}}
\newcommand{\G}{{\mathcal G}}

\newtheorem{theorem}{\bf Theorem}
\newtheorem{proposition}{\bf Proposition}
\newtheorem{lemma}{\bf Lemma}

\newtheorem{definition}{\bf Definition}
\newtheorem{remark}{Remark}

\newtheorem{example}{\bf Example}

\newtheorem{claim}{\bf Claim}

\newtheorem{notation}{\bf Notation}

\begin{document}

\title[Degenerate Foliations]{Foliations with Degenerate Gauss maps on $\mathbb P^4$}

\author{T. FASSARELLA}
\address{Departamento de Matemática -- CCE \\
Universidade Federal do Espírito Santo \\
Av.Fernando Ferrari 514 -- Vitória\\
29075-910 ES Brasil} \email{thiago@impa.br}

 \subjclass{}

\begin{abstract}
We obtain a classification of codimension one holomorphic foliations on $\mathbb P^4$ with degenerate Gauss maps. 
\end{abstract}

\maketitle

\section{Introduction}
Varieties  with degenerate Gauss maps in complex projective spaces and their differential--geometric properties has already been considered by Segre in \cite{S,S2}.  That study was renewed by Griffiths--Harris in \cite{GH} and concepts like the second fundamental form which can be seen as the differential of the Gauss map are discussed in modern language. It presents a growing literature. See \cite{AG,FP,IL} for example.

In the current paper we study varieties with degenerate Gauss maps when are posed in a family, more precisely, we deal with codimension one holomorphic foliations on the complex projective space $\mathbb P^n$ which every leaf has degenerate Gauss map. It happens when the global Gauss map associated to a foliations $\f$ on $\mathbb P^n$ is degenerate. In this case we shall say $\f$ is degenerate.

It is not difficult to see that every curve with degenerate Gauss map is a line and then any degenerate foliation on $\mathbb P^2$ is a pencil of lines. A surface in $\mathbb P^3$ with degenerate Gauss map is either a plane, a cone, or the tangential surface of a curve. According to Cerveau-Lins Neto \cite{CLn} there are two non-exclusive
possibilities for a degenerate foliation on $\mathbb P^3$:
\begin{enumerate}
\item $\f$ has a rational first integral;
\item $\f$ is a linear pull-back of some foliation on $\mathbb P^2$.
\end{enumerate}
Moreover, if $\f$ is not a linear pull-back of some foliation on $\mathbb P^2$, then the rational first integral can be given explicitly and the leaves of $\f$ are cones over a fixed curve. If $\f$ is a linear pull-back of some foliation on $\mathbb P^2$, then the leaves are cones over the leaves of such foliation on $\mathbb P^2$.

Such classification of degenerate foliations on $\mathbb P^3$  was used by Cerveau--Lins Neto in order to understand the irreducible components of the space $\f(n,d)$ of codimension one holomorphic foliations on $\mathbb P^n$ with fixed degree $d=2$. Their analysis relies on one hand in the study of degenerate foliations  on $\mathbb P^3$ and on the other hand in the Dulac's classification of the foliations with degree 2 on $\mathbb P^2$ which has a Morse type singularity. This shows that a relevant role in the understanding of the irreducible components of the space of codimension one holomorphic foliations is played by the analysis of the foliations with degenerate Gauss maps. 

In this work we will obtain a similar classification as in $\mathbb P^3$, but considerably more complicate, of the foliations with degenerate Gauss maps on $\mathbb P^4$. We will see that,  in contrast with lower dimension, there are examples in $\mathbb P^4$ which neither  have rational first integral nor are linear pull-back of some foliation on $\mathbb P^2$. The examples that arise in our classification are separated in three classes according to whether the leaves of $\f$ are cones, joins or bands.  Some of these examples have an interesting structure which are different from the previously known, due to the fact that the degenerate foliation is completely determined by the foliation defined by the fibers of the Gauss map.

The paper is organized as the following: In \S \ref{S:basic} we recall the basic definition of codimension $q$ holomorphic foliations on $\mathbb P^n$ and their Gauss maps. In \S \ref{S:lin} we show the linearity of the fibers of the Gauss map. In \S \ref{S:kplanes} we recall some basic property concerning foliations by $k$-planes. In \S \ref{S:Degenerate}  we prove the Theorem \ref{T:focais} which relates the focal points of the foliation determined by the fibers of the Gauss map with the focal points of a leaf of $\f$ and use this to give a shorter proof of the classification of degenerate foliations on $\mathbb P^3$. In \S \ref{S:degenerate4} we give the examples of degenerate foliations on $\mathbb P^4$ - which possibility do not have rational first integral - that will arise in our classification and prove the Theorem \ref{T:teor_classif}, our main result.

\bigskip

\noindent{\textbf{Acknowledgement:}} 
My warmest thanks to Prof. Jorge Vitório Pereira for his substantial help in the elaboration of this paper. I thank also to IMPA where a  great part of this paper was written.

\section{Basic Definitions}\label{S:basic}
A \textit {codimension $q$ singular holomorphic foliation on $\mathbb P^n$}, from now on,
just codimension $q$ foliation $\f$ on $\mathbb P^n$ is
determined by a line bundle $\mathcal L$ and an element $\omega \in
\mathrm H^0(\mathbb P^n, \Omega^q_{\mathbb P^n} \otimes \mathcal L)$ satisfying
\begin{enumerate}
\item[(i)] $\mathrm{codim}(\sing(\omega)) \ge 2$ where $\sing(\omega)
= \{ x \in \mathbb P^n \, \vert \, \omega(x) = 0 \}$;
\item[(ii)] $\omega$ is integrable.
\end{enumerate}
By definition $\omega$ is integrable if and only if for every point $p\in \mathbb P^n\backslash \sing(\omega)$ there exist a neighborhood $V\subset \mathbb P^n$ of $p$ and  $1$-forms $\alpha_1$,...,$\alpha_q \in \Omega^1(V)$ such that
\begin{eqnarray*}
\omega|_V=\alpha_1\wedge \ldots \wedge\alpha_q \quad \text{and} \quad d\alpha_i \wedge \omega|_V=0 \quad \forall\; i=1,...,q.
\end{eqnarray*}

The \textit{singular set} of $\f$, for short $\sing(\f)$, is by definition equal to  $\sing(\omega)$. The integrability
condition (ii) determines in an analytic neighborhood of every regular point $p$, i.e.
$p \in \mathbb P^n \setminus \sing(\f)$ a holomorphic fibration with
relative tangent sheaf coinciding with the subsheaf of $T\mathbb P^n$
determined by the kernel of $\omega$. Analytic continuation of the
fibers of this fibration describes the \textit{leaves} of $\f$.

For a \textit{codimension one foliation} $\f$ on $\mathbb P^n$ (i.e. $q=1)$ the integrability (ii) is equivalent to
\[
\omega\wedge d\omega=0.
\] 

The \textit{Gauss map} of a codimension $q$ foliation $\f$ on $\mathbb P^n$ is the rational map
\begin{eqnarray*}
\GF : \mathbb P^n &\dashrightarrow& \mathbb G(n-q,n)\,  \\
p &\mapsto& T_p \f
\end{eqnarray*}
where $T_p \f$ is the projective tangent space of the leaf
of $\f$ through $p$ and $\mathbb G(n-q,n)$ is the  Grassmannian of $(n-q)$-planes in $\mathbb P^n$. In this paper we will focus on the  Gauss map of codimension one foliations. Codimension $q$ foliations, $2\le q \le n-1$, will appear only in the definition of foliations by $k$-planes in \S \ref{S:kplanes}. 

The Gauss map for codimension one foliations was also considered in \cite{PY} to give upper bounds for the dimension of certain resonance varieties and in \cite{FPe} for the study of the degree of polar transformations.

\section{Linearity of the Fibers of the Gauss Map}\label{S:lin}
The Gauss map of a codimension one foliation $\f$ on $\mathbb P^n$ is the rational map
\begin{eqnarray*}
\GF : \mathbb P^n &\dashrightarrow& \check
{\mathbb P}^n \,  \\
p &\mapsto& T_p \f
\end{eqnarray*}
where $T_p \f$ is the projective tangent space of the leaf
of $\f$ through $p$.

The degree of a codimension one foliation $\f$ on $\mathbb P^n$, for short ${\rm{deg}}(\f)$, is geometrically defined as the number of tangencies of $\mathcal F$ with a generic line $\ell
\subset \mathbb P^n$. If $\iota: \ell \to \mathbb P^n$ is the
inclusion of such a line then the degree of $\mathcal F$ is the
degree of the zero divisor of the twisted $1$-form  $\iota^*\omega
\in \mathrm H^0(\mathbb \ell, \Omega^1_{\ell} \otimes \mathcal L_{|
\ell})$. Thus the degree of $\mathcal F$ is nothing more than $
\deg(\mathcal L) -2$. It follows from  Euler's sequence that a $1$-form $\omega \in
\mathrm H^0(\mathbb P^n ,\Omega^1 ( \deg(\mathcal F) + 2) )$ can be
interpreted as  a homogeneous $1$-form on $\mathbb C^{n+1}$, still
denoted by $\omega$,
\[
\omega = \sum_{i=0}^n A_i dx_i
\]
with the coefficients $A_i$ being homogenous polynomials of degree
$\deg(\mathcal F) + 1$ and satisfying Euler's relation $ i_R \omega
= 0,$ where  $i_R$ stands for the interior product with the radial
(or Euler's) vector field $R = \sum_{i=0}^n x_i
\frac{\partial}{\partial x_i}$.

If we interpret $[dx_0: \ldots : dx_n]$ as projective coordinates of
$\check{\mathbb P}^n$, then the Gauss map of the corresponding
 $\mathcal F$  can be written  as   $\GF(p)=[A_0(p): \ldots: A_n(p)]$.

If $p$ is a regular point of $\f$ and $M$ is a germ of  leaf of $\f$ through $p$, i.e. a germ of codimension one complex variety in which is invariant by $\f$. We consider the Gauss map associated to $M$
\begin{eqnarray*}
G_M: M &\longrightarrow& \check{\mathbb P}^n \\
         p &\mapsto& T_pM.
\end{eqnarray*}
It is clear that $\GF$ restricted to $M$ coincides with $G_M$. The following lemma relates their ranks.

\begin{lemma}\label{L:fiber}
Let $\f$ be a codimension one foliation on $\mathbb P^n$ and $M$ be a germ of leaf of $\f$  through a regular point $p$. Then
\begin{eqnarray*}
\rank(d\GF(p))=\rank(dG_M(p)) + 1.
\end{eqnarray*}
\end{lemma}
\begin{proof}
We consider 
\[
\omega = \sum_{i=0}^n A_i dx_i,
\]
a homogeneous 1-form on $\mathbb C^{n+1}$ representing $\f$. 

After a linear change of variables one can suppose in the affine coordinate system $(z_1,...,z_n)=(\frac{x_0}{x_n},...,\frac{x_{n-1}}{x_n}) \in \mathbb C^n\subset\mathbb P^n$ that 
\[
p=(0,...,0) \,\, \text{and} \,\, T_p\f=\{z_n=0\}. 
\]

In this coordinate system $q=(z_1,...,z_n)$, the Gauss map is given by
\begin{eqnarray*}
\GF(q)=\frac{1}{a_{n}(q)}\left(a_1(q),...,a_{n-1}(q),-\sum_{i=1}^{n}z_ia_i(q)\right),
\end{eqnarray*}
where $a_i(q)=A_{i-1}(z_1,...,z_n,1)$, $i=1,...,n$.

Therefore, it follows
\[
d\GF(0)\left(\frac{\partial}{\partial z_n}\right) \notin d\GF(0)(T_0\f).
\]
Since $d\GF(0)|_{T_0\f}=dG_M(0)$, this concludes the proof. 
\end{proof}

Suppose $\f$ has degenerate Gauss map, i.e. the differential $d\GF(p)$  has constant rank $k$ for some $1\le k <n$ outside a proper algebraic subset $S\subset \mathbb P^n$. In this case we shall say $\GF$ has rank $k$. By the Implicit Function Theorem there is a regular holomorphic foliation $\mathcal G$ on $\mathbb P^n\setminus S$ so that each leaf of $\mathcal G$ has dimension $n-k$ and $\GF$ is constant on each leaf. Such foliation $\mathcal G$ extend to a codimension $k$ foliation on $\mathbb P^n$ in which we still denote by $\mathcal G$. We shall refer to $\mathcal G$ as the \textit{foliation determined by the fibers of $\GF$}. The following proposition shows that the leaves of $\mathcal G$ are linear spaces.

\begin{proposition}\label{P:fibra_linear_fol}
Let $\f$ be a codimension one foliation on $\mathbb P^n$, $n\ge 2$. If $\GF$ has rank k, $1\le k < n$, then a generic fiber of $\GF$ is a union of open subsets of linear spaces of dimension $n-k$.   
\end{proposition}

\begin{proof}
Let $M$ be the germ of a leaf of $\f$. By  Lemma \ref{L:fiber}, the generic fiber of $\GF$ and $G_M$ has the same dimension. To conclude the proposition one has just to invoke the classical results on the Gauss map $G_M$ (cf. \cite{S} and \cite{GH} for a local approach).
\end{proof}

\section{Foliations by \textit{k}-Planes}\label{S:kplanes}

A \textit{foliation by $k$-planes} $\mathcal G$ on $\mathbb P^n$, $1\le k \le n-1$, by definition is a codimension $(n-k)$ foliation on $\mathbb P^n$ which every leaf is contained in a linear space of dimension $k$. The algebraic closure $B'\subset \mathbb G(k,n)$ of the set of invariant $k$-planes by $\mathcal G$  defines a subvariety of dimension $n-k$. Reciprocally, a subvariety $B'\subset  \mathbb G(k,n)$ which for a generic point $p \in \mathbb P^n$ there is exactly one element of $B'$ passing through $p$, defines a foliation by $k$-planes $\mathcal G$ on $\mathbb P^n$. We shall say $\mathcal G$ \textit{is determined by} $B'$.

We consider the variety
\[
\Lambda '=\{(p,L) \in \mathbb P^n \times B '\,|\, p\in L\}.
\]
together with the natural projections $g ':\Lambda '\longrightarrow B '$ and $f ':\Lambda '\longrightarrow \mathbb P^n$. Then we have the following diagram
\[
 \xymatrix{
\Lambda:=B \times_{B '}\ar@/^0.5cm/@{->}[rr]^{f}\Lambda ' \ar[d]_{g} \ar[r]_{\phi^*}  & \Lambda '  \ar[d]_{g '} \ar[r]_{f '}  & \mathbb P^n
  \\
  B \ar[r]_{\phi } & B '}
\]
where $\phi$ is a resolution of $B'$ and $f=f '\circ \phi^*$. Notice that $\Lambda$ is a vector bundle over the nonsingular variety $B$, in particular it is a smooth variety.

We denote by $R_{\mathcal G}$ the divisor of ramification points of $f$. We can write 
\[
R_{\G}=H_{\G}+V_{\G}
\]
where $g_{|_{|H_{\G}|}}$ has rank $n-k$ in a dense open subset and $g_{|_{|V_{\G}|}}$ is not dominant (where $|\;\;|$ denotes the support of the divisor). Notice that the irreducible components of the support of $H_{\G}$ can be seen as \textit{horizontal components} (and $V_{\G}$ as vertical components) with respect to the $k$-planes fibers of $g$. Since $f$ is a parameterization of the leaves of $\G$ outside $f(|R_{\G}|)$, it follows $f(|R_{\G}|)=\sing(\G)$.

\begin{definition}\rm
If $E$ is an irreducible component of the support of $H_{\G}$, then $f(E)$ is called \textit{fundamental component of} $\G$. The union of all fundamental components is called \textit{fundamental set} of $\G$ and it is denoted by $\Delta_{\G}$.
\end{definition}

\begin{figure}[h]
\centering
\includegraphics[width=0.4\textwidth]{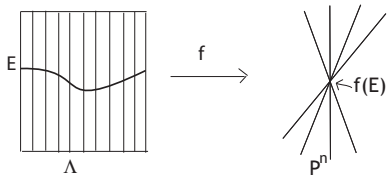}
\caption{Fundamental component.}
\end{figure}

If $L\in B$, we still denote by $L \subset \mathbb P^n$ the respective $k$-plane interpreted as subset of $\mathbb P^n$. If $L$ is an invariant $k$-plane by $\G$  is not contained in the singular set of $\G$, that is, $g^{-1}(L)$ is not contained in the support of $V_{\G}$, then the restriction of $H_{\G}$ to $g^{-1}(L)$ defines a hypersurface in $g^{-1}(L)$ which we denote by $\Delta_{\G}^*(L)$. Notice that 
\[
f|_{g^{-1}(L)}:g^{-1}(L) \longrightarrow L
\]
is an isomorphism. Thus the direct image $f_*(\Delta_{\G}^*(L))$ is a divisor in $L \subset \mathbb P^n$ which we denote by $\Delta_{\G}(L)$. We shall refer to $\Delta_{\G}(L)$ as the \textit{focal points of $\G$ in L}. We will proof in \S \ref{Sub:focalpoints} that $\Delta_{\G}(L)$ has degree $n-k$. That result can also be found in \cite{Dep3}, it is a classical known fact about congruences.

\medskip

The subvariety $B'\subset \mathbb G(k,n)$ which defines a foliation by $k$-planes is classically known as \textit{congruence of $k$-planes  on $\mathbb P^n$ of order one}. The congruence of lines on  $\mathbb P^n$ (i.e. $k=1$) was already considered by Kummer in \cite{Kum} where a classification was obtained when $n=3$. We state this classification in \S \ref{Su:lines3} and use it to obtain a shorter proof of Cerveau--Lins Neto's  classification of degenerate foliations on $\mathbb P^3$. The congruence of lines on $\mathbb P^4$ was studied by Marletta in \cite{Mar,Mar2} and also by De Poi in \cite{Dep2,Dep3,Dep4} seeking to complete Marletta's classification. We will not use such classification for the study of degenerate foliations on $\mathbb P^4$ in \S \ref{S:degenerate4}, just some general results. 

\subsection{Foliations by Lines}\label{Sub:foliationlines}
Now, let us recall some known facts (see \cite{Dep3} for example) about foliations by lines (i.e. $k=1$) that will be useful later.

\begin{lemma}\label{L:folhe_retas_isol}
Let $\G$ be a foliation by lines on $\mathbb P^n$. If $\G$ has an isolate singularity $p$, then it is the radial foliation from $p$, that is,  the foliation determined by the lines passing through  $p$. 
\end{lemma}

\begin{proof}
Let $E$ be an irreducible component of the support of $R_{\G}$ such that $p\in f(E)$. Since $p$ is isolate we conclude that $f(E)=\{p\}$.

If $E$ is contained in the support of $V_{\G}$, then $f(E)\subset\sing(\G)$ contains at least one line $f(g^{-1}(L))$.  Since $f_{|{g^{-1}(L)}}$ is an isomorphism,  $p$ would not be an isolate singularity. 

Therefore we may assume that $E$ is contained in the support of $H_{\G}$.  It follows that $B '\subset \sigma$, where $B'\subset \mathbb G(1,n)$ denote the closure of the family of invariant lines by $\G$ and $\sigma \subset \mathbb G(1,n)$ denote the family of lines passing through $p$. Since 
\[
{\rm{dim}}(B ')={\rm{dim}}(\sigma)=n-1, 
\]
the result follows.    
\end{proof}

\begin{lemma}\label{L:det}
Let $\G$ be a foliation by lines on $\mathbb P^n$. Let us  write 
\begin{eqnarray*}
H_{\G}=\sum m(E)\cdotp E,
\end{eqnarray*}
where $E$ is an irreducible hypersurface in $\Lambda$. We claim that 
\[
m(E)\ge n- {\rm{dim}}f(E)-1.
\]
\end{lemma}

\begin{proof}
We take $q\in E$ a smooth point such that $p=f(q)$ is a smooth point of $f(E)$. Up to a local biholomorphism we may assume that $f$ can be written locally as
\[
f=(a_1,...,a_{d_E},x_na_{d_E+1},...,x_na_n)
\]
where $a_i \in \mathcal O_q$, for all $i=1,...,n$, $x_n\in \mathcal O_q$ is the local equation of $E$ and $d_E = {\rm{dim}} f(E)$. It follows 
\[
{\rm{det}}Jf=x_n^{n-d_E-1}\cdot Q,
\]
where $Q \in \mathcal O_q$. Since such biholomorphism induces an isomorphism in the local rings, this concludes the lemma.
\end{proof}

\begin{definition}\rm\label{D:projalong}
Let $\G$ be a foliation by lines on $\mathbb P^n$ and $E\subset \Lambda$  an irreducible component of the support of $H_{\G}$.
\begin{enumerate}
\item If \, $g^{-1}(L)\cap E$ is a set of single point for generic $L \in B$, we may consider the rational map 
\begin{eqnarray*}
\phi_E := f\circ \tau_E \circ f^{-1}:\mathbb P^n \dashrightarrow f(E),
\end{eqnarray*}
where $\tau_E:\Lambda \longrightarrow E$ is the natural projection along the lines $g^{-1}(L)\subset \Lambda$. We shall refer to $\phi_E$ as the \textit{projection along the lines of $\G$}.

\item The fibers of the restriction $f_{|E}:E \longrightarrow f(E)$ define a foliation $\eta_E$ on $E$ of dimension $n-d_E-1$. We shall refer to $\eta_E$ as the \textit{foliation defined by $f$ on $E$}.
\end{enumerate}
\end{definition}

\begin{proposition}\label{P:depoi}
Let $\G$  be a foliation by lines on $\mathbb P^n$, $n=3,4$. If $E$ is an irreducible component of the support of $H_{\G}$ such that ${\rm{dim}}f(E)\ge 1$, then 
\begin{eqnarray*}
{\rm{ker}}df(x) \subset T_xE, 
\end{eqnarray*}
for a generic point $x\in E$.
\end{proposition}
\begin{proof}
Notice that if $T_xE$ does not contain ${\rm{ker}}df(x)$, then
\begin{eqnarray*}
T_x\Lambda=T_xE + {\rm{ker}}df(x).
\end{eqnarray*}
Therefore $df(x)(T_x\Lambda)\subset T_{f(x)}f(E)$, i.e. the line $f(g^{-1}(g(x)))$ is tangent to $f(E)$ in $f(x)$. If ${\rm{dim}}f(E)= 1$, then the closure $B '\subset \mathbb G(1,n)$ of invariant lines  by $\G$  is a family of tangent lines to a curve. But this family does not fill up the space. If ${\rm{dim}}f(E)= 2$, then  $B '$ is the family of tangent lines to a surface. But it cannot define a foliation by lines according to \cite{Dep2}, Proposition $1.3$. 
\end{proof}

We need introduce more notations:
\begin{notation}\label{N:esp_linear}\rm$\empty$

\begin{enumerate}
\item If $F\subset\mathbb P^n$ is a linear space, then $F^{\lor}$ denote its dual, i.e. the projective space associated to the set of hyperplanes in $F$.
\item If $E\subset F$ is a linear subspace, then $E_F^*$ denote the set of hyperplanes $\pi\in F^{\lor}$ containing $E$. If $F=\mathbb P^n$ we just denote $E^*$.
\end{enumerate}
\end{notation}

\subsection{Foliations by Lines on $\mathbb P^3$}\label{Su:lines3}
The foliation by lines were considered by Kummer in \cite{Kum}, where a classification in $\mathbb P^3$ was obtained. Alternative proofs of the Kummer's classification has been given by many other authors (cf. \cite{Ran,ZILO, Dep, C}). 

We state such classification here:

\begin{theorem}\label{T:retas3}(Kummer, \cite{Kum}).
Let $\G$ be a foliation by lines on $\mathbb P^3$. We have one of the following possibilities:

\begin{enumerate}
\item The foliation $\G$ is determined by the family of lines passing through a point $p \in \mathbb P^3$. 
\item The foliation $\G$ is determined by the family of secant lines to a twisted cubic.
\item There is a line $L$ and a non-constant morphism $\psi:L^* \longrightarrow L$ such that $\G$ is determined by the family of lines
\begin{eqnarray*}
\overline{\cup_{\pi\in L^*}\psi(\pi)_{\pi}^*} \;\;.
\end{eqnarray*}
\item There is a line $L$ and a rational curve $C$  such that for generic $\pi\in L^*$, then $\left( \pi \cap C\right) \backslash L $ is a set of a single point and $\G$ is determined by the family of lines which intersect $L$ and $C$.
\end{enumerate}
\end{theorem}

\subsection{Foliations by Planes on $\mathbb P^4$}\label{Su:planes4} 
The classification of the surfaces of order one in $\mathbb G(k,n)$, that is, of the family of $k$-planes in $\mathbb P^n$ which a generic $(n-k-2)$-plane meets only one $k$-plane of the family was obtained by Z.Ran in \cite{Ran}. In particular, he gave a classification of the foliation by planes $(i.e.\; k=2)$ on $\mathbb P^4$ which we state here. This classification will be useful later.

\begin{theorem}\label{T:class_folh_planos}(Ran, \cite{Ran}).
Let $\G$ be a foliation by planes on $\mathbb P^4$. We have one of the following possibilities:
\begin{enumerate}
\item The foliation $\G$ is determined by the family of planes which contains a line.
\item There is a plane $\Sigma$ and a morphism $\psi:\Sigma^* \longrightarrow \mathbb G(1,4)$ such that the line $\psi(\pi)$ is contained in $\pi$ for all hyperplane $\pi\in\Sigma^*$ and $\G$ is determined by the family of planes
\begin{eqnarray*}
\overline{\cup_{\pi\in\Sigma^*}\cup_{L\in(\psi(\pi))_{\pi}^*}L}.
\end{eqnarray*}
\item There is a cubic irreducible surface $S$ and a family of dimension two of conics (distinct of double lines) in $S$ such that $\G$ is determined by the family of planes containing such conics.
\end{enumerate}
\end{theorem}

\section{Degenerate Foliations}\label{S:Degenerate}
Given a codimension one foliation $\f$ on $\mathbb P^n$ one can naturally associate to it the Gauss map $\GF$. When $\GF$ has rank $n-k$, $1 \le k < n$, it follows from Proposition \ref{P:fibra_linear_fol} that its fibers define a foliation $\G$ by $k$-planes  on $\mathbb P^n$.

\begin{definition}\rm
A codimension one foliation on $\mathbb P^n$ is called \textit{degenerate} if it has degenerate Gauss map. The \textit{rank} of $\f$ is defined as the rank of $\GF$.
\end{definition}

\subsection{Focal Points}\label{Sub:focalpoints}
Let $\f$ be a codimension one foliation of rank $n-k$, $1\le k < n$. Let $p\in\mathbb P^n$ be a regular point of $\f$, $U\subset \mathbb P^n$  a small neighbourhood of $p$ in which does not intersect $\sing(\f)$ and  $M\subset U$ a codimension one complex variety which is invariant by $\f$. We consider  
\[
M_0:=\cup_{p\in M}L_p ,
\]
where $L_p$ is the invariant $k$-plane by $\G$ passing through $p$. We shall refer to $M_0$ as the \textit{saturation of $M$ by $\G$}. To obtain a definition of the focal points for $M_0$ we take the maps $f_{M_0}$ and $g_{M_0}$ defined in the following diagram
\[
 \xymatrix{
\Lambda_{M_0}\ar@/^0.5cm/@{->}[rr]^{f_{M_0}:=f\circ i_1}   \ar[d]_{g_{M_0}:=g\circ i_1} \ar[r]_{i_1}  & \Lambda   \ar[d]_{g} \ar[r]_{f}  & \mathbb P^n
  \\
  B_{M_0} \ar[r]_{i_2} & B }
\]
where $\Lambda_{M_0}=f^{-1}(M_0)$, $B_{M_0}=g(\Lambda_{M_0})$  and $i_1$, $i_2$ are natural inclusions. 

We denote by $R_{M_0}$ the divisor of ramification points of $f_{M_0}$. We can write 
\[
R_{M_0}=H_{M_0}+V_{M_0}
\]
where $(g_{M_0}){|_{{|H_{M_0}|}}}$ has rank $n-k-1$ in a dense open subset and $(g_{M_0}){|_{|V_{M_0}|}}$ is not dominant.

If $L$ is an invariant $k$-plane by $\G$ is not contained in $\sing(\G)$, that is, $g_{M_0}^{-1}(L)$ is not contained in the support of $R_{M_0}$, then the restriction of $H_{M_0}$ to $g_{M_0}^{-1}(L)$ defines a divisor in $g_{M_0}^{-1}(L)$ which we denote by $\Delta_{M_0}^*(L)$. Notice that 
\[
f_{M_0}|_{g_{M_0}^{-1}(L)}:g_{M_0}^{-1}(L) \longrightarrow L
\]
is an isomorphism. Then the direct image 
\[
(f_{M_0})_*(\Delta_{M_0}^*(L))
\]
is a divisor in $L \subset \mathbb P^n$ which we denote by $\Delta_{M_0}(L)$. We shall refer to $\Delta_{M_0}(L)$ as the \textit{focal points of $M_0$ in L}. 

The following theorem relates the focal points of $\G$ and the focal points of $M_0$. In order to prove it, we use the same idea as in \cite[\S 2.2.4]{FP} to look at the focal points.

\begin{theorem}\label{T:focais}
With notations as above, if $L \subset \mathbb P^n$ is an invariant $k$-plane  by $\G$ is not contained in $\sing(\G)$, then we have the following statements:
\begin{enumerate}
\item $\Delta_{M_0}(L)$ is a divisor in $L$ of degree $n-k-1$;
\item $\Delta_{\G}(L)$ is a divisor in $L$ of degree $n-k$;
\item $\Delta_{M_0}(L)\le\Delta_{\G}(L)$.
\end{enumerate} 
\end{theorem}

\begin{proof}
Since $f|_{g^{-1}(L)}:g^{-1}(L) \longrightarrow L$ is an isomorphism we can suppose that
\begin{eqnarray*}
f^*\Delta_{\G}(L)={\rm{V}}(F_{\G}(L)) \;\; \text{and} \;\;
f^*\Delta_{M_0}(L)={\rm{V}}(F_{M_0}(L)),
\end{eqnarray*}
where $F_{\G}(L)$, $F_{M_0}(L)$ are homogeneous polynomial in $k+1$ variables.

A  local parameterization $\psi:\mathbb D^{n-k} \longrightarrow B$  of a neighbourhood of $L$ can be explicitly written as
$\psi(s)=\mathbb P({\rm{span}}\{ \alpha_0(s),...,\alpha_k(s)\})$,
where $\alpha_0,..., \alpha_k : \mathbb D^{n-k} \longrightarrow \mathbb C^{n+1}$ are holomorphic functions. Composing $f$ with the following local parameterization 
\begin{eqnarray*}
\lambda:\mathbb D^{n-k}\times\mathbb P^k &\longrightarrow& \Lambda\\
                 (s,t)    &\mapsto&   ([t_0\alpha_0(s)+...+t_k\alpha_k(s)],\psi(s)),  \;\; t=(t_0:...:t_k),
\end{eqnarray*}
we may assume that $f$ can be explicitly written as
\begin{eqnarray*}
f:\mathbb D^{n-k}\times\mathbb P^k &\longrightarrow& \mathbb P^n\\
 (s,t)    &\mapsto&   [t_0\alpha_0(s)+...+t_k\alpha_k(s)],  \;\; t=(t_0:...:t_k).
\end{eqnarray*}
Since $L$ is not contained in $\sing(\G)$ one can suppose that
\begin{enumerate}
\item $f$ is an isomorphism in $(s,e_0)$ (where $e_0=(1:0:...:0)$), for all  $s\in \mathbb D^{n-k}$ and we take $\Sigma:=f(\mathbb D^{n-k}\times\{e_0\})$ transversal section to $\G$; 
\item $f(\mathbb D^{n-k-1}\times\{s_{n-k}\}\times\{e_0\})$ is invariant by $\f|_{\Sigma}$ for all $s_{n-k} \in \mathbb D$.
\end{enumerate} 
Therefore, it follows from the above property $(2)$ that the restrictions
\begin{eqnarray*}
f_{s_{n-k}}:\mathbb D^{n-k-1}\times\mathbb P^k &\longrightarrow& \mathbb P^n \\
 (\overline{s},t)  &\mapsto&  f(s,t) , \;\; s=(\overline{s},s_{n-k}),
\end{eqnarray*} 
are local parametrizations of the leaves of $\f$. Then by definition, $f^*\Delta_{\G}(L)$ is given by the ramification points of $f$ and $f^*\Delta_{M_0}(L)$ is given by the ramification points of $f_{s_{n-k}}$.

Let us denote 
\[
V:={\rm{span}}\left\{\alpha_0(s),...,\alpha_k(s),\sum_{i=0}^{k}t_i\frac{\partial\alpha_i}{\partial s_1}(s),...,\sum_{i=0}^{k}t_i\frac{\partial\alpha_i}{\partial s_{n-k}}(s)\right\}.
\]
We notice that ${\rm{Im}}\,df_{(s,t)}=\mathbb P V$. Hence property $(1)$ above implies that $f$ is an isomorphism in $(s,t)$ if and only if the equality in the following equation is hold
\begin{eqnarray*}
V\subset {\rm{span}}\left\{\alpha_0(s),...,\alpha_k(s),\frac{\partial\alpha_0}{\partial s_1}(s),...,\frac{\partial\alpha_0}{\partial s_{n-k}}(s)\right\} =\mathbb C^{n+1}.
\end{eqnarray*}

If we denote
\begin{eqnarray*}
\frac{\partial\alpha_i}{\partial s}(s):=\left(\frac{\partial\alpha_i}{\partial s_1}(s),...,\frac{\partial\alpha_i}{\partial s_{n-k}}(s)\right),\;\;i=0,...,k, 
\end{eqnarray*}    
then, there exists matrix $A_i(s) \in M((n-k);\mathbb C)$, $i=0,...,k$, $(A_0(s)=Id_{n-k})$ such that 
\begin{eqnarray*}
\frac{\partial\alpha_i}{\partial s}(s)=A_i(s)\frac{\partial\alpha_0}{\partial s}(s) \;\;\;\;\;{\rm{mod}} \;{\rm{span}}\{ \alpha_0(s),...,\alpha_k(s)\}.
\end{eqnarray*}
Therefore $f$ is an isomorphism in $(s,t)$ if and only if
\begin{eqnarray*}
{\rm{det}}(t_0A_0(s)+...+t_k A_k(s))=0.
\end{eqnarray*}
This implies that $F_{\G}(L)={\rm{det}}(t_0A_0(s)+...+t_kA_k(s))$ and hence is a homogeneous polynomials in $t_0,...,t_k$ of degree $n-k$. This proves the statement $(2)$ of the theorem.

On the other hand, the following equation  follows from the hypothesis that the tangent space of $\f$ is constant along the leaves of $\G$ 
\begin{eqnarray*}
&{\rm{span}}\left\{\alpha_0(s),...,\alpha_k(s),\sum_{i=0}^{k}t_i\frac{\partial\alpha_i}{\partial s_1}(s),...,\sum_{i=0}^{k}t_i\frac{\partial\alpha_i}{\partial s_{n-k-1}}(s)\right\}& \subset \\ &{\rm{span}}\left\{\alpha_0(s),...,\alpha_{k}(s),
\frac{\partial\alpha_0}{\partial s_1}(s),...,\frac{\partial\alpha_0}{\partial s_{n-k-1}}(s)\right\}.&
\end{eqnarray*}
We have equality if and only if $f_{s_{n-1}}$ is an immersion in $(\overline{s},t)$. Therefore there exists matrix $\overline{A}_i(\overline{s}) \in M((n-k-1);\mathbb C)$, $i=0,...,k$, ($\overline{A}_0(\overline{s})=Id_{n-k-1}$) such that 
\begin{eqnarray*}
\frac{\partial\alpha_i}{\partial \overline{s}}(s)=\overline{A}_i\frac{\partial\alpha_0}{\partial \overline{s}}(s) \;\;\;\;\;{\rm{mod}} \;{\rm{span}}\{ \alpha_0(s),...,\alpha_k(s)\}.
\end{eqnarray*}
Hence $f_{s_{n-1}}$ is not a immersion in $(\overline{s},t)$ if and only if
\begin{eqnarray*}
{\rm{det}}(t_0\overline{A}_0(\overline{s})+...+t_k\overline{A}_k(\overline{s}))=0.
\end{eqnarray*}
This implies that  $F_{M_0}(L)={\rm{det}}(t_0\overline{A}_0(\overline{s})+...+t_k\overline{A}_1(\overline{s}))$ and hence is a homogeneous polynomials in $t_0,...,t_k$ of degree $n-k-1$. This concludes the proof of the statement $(1)$. Since
\[
t_0A_0(s)+...+t_kA_k(s)= \left( \begin{array}{c|c} t_0\overline{A}_0(\overline{s})+...+t_k\overline{A}_k(\overline{s}) & 0_{(n-k-1)\times 1}  \\ \hline a_{n-k1}(s)\dots a_{n-kn-k-1}(s) & a_{n-kn-k} (s) \end{array} \right) .
\]
for some holomorphic functions $a_{ij}:\mathbb D^{n-1}\longrightarrow
 \mathbb C$, it follows $F_{M_0}(L)$ divides $F_{\G}(L)$. This concludes the proof of the theorem.
\end{proof}

\begin{remark}\rm $\empty$
\begin{enumerate}
\item[(a)] The statement $(1)$ of the Theorem \ref{T:focais} is a known fact about varieties with degenerate Gauss maps (see \cite[Proposition 2.2.4]{FP} for example). 
\item[(b)] In the case $k=1$, the statement $(2)$ of the Theorem \ref{T:focais} is a known fact about congruence of lines (see \cite[Theorem 5]{Dep3} for example).
\end{enumerate}
\end{remark}

\subsection{Degenerate Foliations on $\mathbb P^3$}\label{Su:degenerate3}
The classification of degenerate foliations on $\mathbb P^3$ plays a key role in the study of the irreducible components of the space of holomorphic foliations of degree two on $\mathbb P^n$ (cf. \cite{CLn}). We give a shorter proof of such classification by using the classification of the foliations by lines on $\mathbb P^3$. The proof contains some main ideas for the classification on $\mathbb P^4$.

\begin{theorem}(Cerveau--Lins Neto, \cite{CLn}).\label{T:class_folh_deg_3}
Let $\f$ be a degenerate foliation on $\mathbb P^3$. We have one of the following possibilities:
\begin{enumerate}
\item $\f$ has a rational first integral;
\item $\f$ is a linear pull-back of some foliation on $\mathbb P^2$.
\end{enumerate}
\end{theorem}
\begin{proof}
If $\f$ has rank one, then the leaves of $\f$  are open subsets of linear spaces (fibers of $\GF$). In this case it is easy to see that  $\f$ is a pencil of hyperplanes. Hence it has a rational first integral.

If $\f$ has rank two, we consider the foliation by lines $\G$ determined by the fibers of $\GF$. If $\G$ has an isolate singularity $p$,  it follows from Lemma \ref{L:folhe_retas_isol} that $\G$ is determined by the family of lines passing through  $p$. Therefore, $\f$ is a linear pull-back of some foliation on $\mathbb P^2$. In fact, let $\psi:\mathbb P^3 \dashrightarrow \mathbb P^2$ be the linear projection from $p$. By hypothesis $\psi$ has fibers which are tangent to $\f$, thus it follows from \cite[Lemma 2.2]{CLnLPT}  (applied to $\psi_{|\left(\mathbb P^3 \backslash \{p\}\right)}$) that there exist a foliation $\eta$ on $\mathbb P^2$ such that $\f=\psi^*(\eta)$. So we are in the case $(2)$ of the statement. From now on we will assume that $\G$ does not have isolate singularities. 

Let $E$ be an irreducible component of $|H_{\G}|$ and $M_0$  be the saturation of a germ of leaf of $\f$ by $\G$. 

On the one hand, the Proposition \ref{P:depoi} implies that (see Definition \ref{D:projalong} for $\eta_E$)
\[
T_p\eta_E={\rm{ker}}(df(p)), \,\,\, \text{for all} \,\,\, p \in E.
\]
On the other hand, since $M_0^*:=f^{-1}(M_0)$ is invariant by $f^*\f$ one has
\[
T_p{\f_E}=T_p(M_0^*\cap E).
\]
By definition of $\Delta_{M_0}^*(L)$ follows
\[
T_p\eta_E=T_p\f_E  \,\,\, \text{if and only if} \,\,\, p \in |\Delta_{M_0}^*(L)|\cap E. 
\]
Therefore $|\Delta_{M_0}^*(L)|\cap E$ is not empty for a saturation $M_0$ of a germ of generic leaf of $\f$ if and only if $\f_E=\eta_E$. It follows from Theorem \ref{T:focais} that $|\Delta_{M_0}^*(L)|\cap E$ is not empty for some $E$, finally one can conclude that $\f_E=\eta_E$ for some irreducible component $E$ of the support of $H_{\G}$.

It follows from Theorem \ref{T:focais} that $\Delta_{\G}(L)$ is a set in $L$ of two points counting multiplicity, so one have three  possibilities:
\begin{enumerate}
\item[(a)] $H_{\G}=2E$; 

\item[(b)] $H_{\G}=E$; 

\item[(c)] $H_{\G}=E_1+E_2$.
\end{enumerate} 

We suppose first $H_{\G}=2E$. By the classification of the foliations by lines on $\mathbb P^3$ we can assume that $f(E)$ is a line and $\G$ is the foliation in the case $(3)$ of the Theorem \ref{T:retas3}. We consider the projection along the lines of $\G$ 
\[
\phi_E:\mathbb P^3 \dashrightarrow f(E).
\]
It follows from the above remark that $\f_{E}=\eta_{E}$,  in particular $\phi_E$ is a rational first integral of $\f$. But the fibers of $\phi_E$ are planes containing $f(E)$, i.e. $\f$ is a pencil of planes containing $f(E)$. Hence it has a rational first integral.

Let now $H_{\G}=E$. It follows from Theorem \ref{T:retas3} that $f(E)$ is a twisted cubic and $\G$ is determined by its secants lines, i.e. $\G$ is the foliation of the case $(2)$ of the Theorem \ref{T:retas3}. Since a generic line $g^{-1}(L)$ intersect $E$ in two points and $\f_E=\eta_E$ it follows that $|\Delta_{M_0}^*(L)|$ is a set of two points. However, this contradicts the fact that a leaf of $\f$ has only one focal point. Therefore this case cannot happen.

It remains to consider the case in which $H_{\G}=E_1+E_2$.  We first prove that $f(E_1)\ne f(E_2)$. In fact, suppose by contradiction that $C:=f(E_1)=f(E_2)$. The cone $F_p$, $p\in C$, determined by the invariant lines by $\G$ passing through $p$ is irreducible because $C$ is an irreducible curve. Let $N_i$ be an irreducible component of ${f|_{E_i}}^{-1}(p)$, $i=1,2$. Then the images by $f$, $f(g^{-1}(g(N_i)))$, $i=1,2$, are distinct irreducible cones determined by invariant lines by $\G$ passing through $p$. This contradicts the fact that $F_p$ is irreducible.

We may assume that $f(E_1)\ne f(E_2)$. By  classification of the foliation by lines on $\mathbb P^3$ we can suppose that $f(E_1)$ is a line and $\G$ is the foliation of the case $(4)$ of the Theorem \ref{T:retas3}. We have that either $\f_{E_1}=\eta_{E_1}$  or  $\f_{E_2}=\eta_{E_2}$. In the both cases,  the projection along  the lines of $\G$
\[
\phi_{E_i}:\mathbb P^3 \dashrightarrow f(E_i), \,\,i=1 \,\,\text{or}\,\,2,
\]
is a rational first integral of $\f$. So that is the case $(1)$ of the statement. It is possible to see that in the first case, the leaves of $\f$ are cones over $f(E_2)$ and in the second case, $\f$ is a pencil of planes containing $f(E_1)$.
\end{proof}

\begin{remark}\rm
If $\f$ is a degenerate foliation on $\mathbb P^3$ in which is not a linear pull-back of some foliation on $\mathbb P^2$, then we are in the case which $H_{\G}=E_1+E_2$, $f(E_1)$ is a line and $\phi_{E_1}$ is a rational first integral of $\f$. It is possible to prove that  $\phi_{E_1}$ is explicitly given by 
\begin{eqnarray*}
\phi_{E_1}:\mathbb P^3 &\dashrightarrow& \mathbb P^1 \\
(x_0:...:x_3) &\mapsto& (x_2P(x_0,x_1)-Q(x_0,x_1):x_3P(x_0,x_1)-R(x_0,x_1)),
\end{eqnarray*}
where $P$, $Q$ and $R$ are homogeneous polynomials in $\mathbb C[x_0,x_1]$, with ${\rm{deg}}(P)+1={\rm{deg}}(Q)={\rm{deg}}(R)$ given a global parameterization of the rational curve $f(E_2)$
\begin{eqnarray*}
\mathbb P^1 &\dashrightarrow& f(E_2) \subset \mathbb P^3 \\
         (x_0:x_1)    &\mapsto& (x_0P(x_0,x_1):x_1P(x_0,x_1):Q(x_0,x_1):R(x_0,x_1)).
\end{eqnarray*}
\end{remark}

\section{Degenerate Foliations on $\mathbb P^4$}\label{S:degenerate4}

If $\f$ is a foliation of rank at most two on $\mathbb P^n$, then we have the following possibilities (cf. \cite{CLn}):
\begin{enumerate}
\item $\f$ has a rational first integral;
\item $\f$ is a linear pull-back of some foliation on $\mathbb P^2$.
\end{enumerate}
In fact, if $\Sigma$ is a generic linear space of dimension three, then $\f|_{\Sigma}$ has also rank at most two. It follows form Theorem \ref{T:class_folh_deg_3} that either $\f|_{\Sigma}$ has a rational first integral or $\f|_{\Sigma}$ is a linear pull-back of some foliation on $\mathbb P^2$. The result follows from \cite[Lemmas $2$ and $3$]{CLn}. Therefore, to obtain a similar classification  as in Theorem \ref{T:class_folh_deg_3} of degenerate foliations on $\mathbb P^4$, it remains to consider foliations of rank three. 

Let $\f$ be a foliation of rank three on $\mathbb P^4$. In this case, the fibers of the Gauss map define a foliation by lines $\G$ on $\mathbb P^4$. Let $p\in\mathbb P^n$ be a regular point of $\f$, $U\subset \mathbb P^n$ be a small neighbourhood of $p$ in which does not intersect $\sing(\f)$ and  $M\subset U$ a codimension one complex variety which is invariant by $\f$. It follows from Lemma \ref{L:fiber} that $M$ is a hypersurface in $\mathbb P^4$ with degenerate Gauss map of rank two. We consider the saturation of $M$ by $\G$
\[
M_0:=\cup_{p\in M}L_p ,
\]
where $L_p$ is the invariant line by $\G$ passing through $p$. It follows from classification of hypersurface in $\mathbb P^4$ of rank two (see \cite{AG} for a local approach and also \cite{MT,Rog} for projective varieties) that one has the following possibilities:
\begin{enumerate}
\item $M_0$ is a cone;
\item $M_0$ is the join of two curves;
\item $M_0$ is an union of  planes (but not osculating) containing the tangent lines to a curve $C$, i.e. $M_0$ is a band over $C$;
\item There exists a curve $C$ and a surface $S$ such that $M_0$ is an union of tangent lines to $S$ intersecting $C$;
\item There exists two surfaces $S_1$,$S_2$ such that $M_0$ is an union of their common tangent lines (possibly $S_1=S_2$).
\end{enumerate}

If $M_0$ is in the cases $(4)$ or $(5)$ for a generic point $p \in \mathbb P^4$, then the invariant lines by $\G$ are tangents to a surface contained in the fundamental set of $\G$. But, it follows from \cite[Proposition 1.3]{Dep2} that it cannot define a foliation by lines on $\mathbb P^4$.  Therefore a generic leaf of $\f$ is either a \textit{cone}, a \textit{join} or a \textit{band}. 

Before giving the examples of degenerate foliations on $\mathbb P^4$, let us give a proposition that will be useful later.

\begin{proposition}\label{P:rankthree}
Let $\f$ be a degenerate foliation on $\mathbb P^4$ which does not have rational first integral. Suppose there is a foliation by lines $\G$ on $\mathbb P^4$ in which is tangent to the foliation determined by the fibers of $\GF$. If $\G$ is not tangent to a foliation by  planes containing a line, then $\f$ has rank three. 
\end{proposition}

\begin{proof}
If $\f$ does not have rank three, then either it has a rational first integral or is a linear pull-back of some foliation on $\mathbb P^2$. By hypothesis one can suppose $\f$ is the pull-back by a linear projection $\psi:\mathbb P^4 \dashrightarrow \mathbb P^2$ with center at the line $L$, of some foliation on $\mathbb P^2$. 

If $\f$ has rank one, then their leaves coincide with the fibers of $\GF$. So, in this case $\f$ is a pencil of hyperplanes containing a line. Hence $\f$ has a rational first integral. Therefore it remains to consider the case which $\f$ has rank two.

Suppose $\f$ has rank two. Since the  foliation determined by the fibers of $\psi$  is tangent to the foliation determined by the fibers of $\GF$, this implies that these two foliations are the same. If $\G$ is a foliation by lines under the hypothesis of the proposition, then $\G$ must be tangent to the fibers of $\psi$. Since the fibers of $\psi$ are planes containing the line $L$, this concludes the proof.
\end{proof}

\subsection{The Examples}\label{Su:theexamples}
In this section we will give the examples of degenerate foliations on $\mathbb P^4$ that possibility does not have rational first integral. Such examples are separated in three classes according the generic leaf of the foliation: cones, joins and bands. In first example the leaves are cones, in second example the leaves are joins and for the rest  the leaves are bands. The examples \ref{E:riccati}, \ref{E:bernoulli} and \ref{E:veronese} are  different from the previously known because the degenerate foliation is completely determined by the foliation by lines defined by the fibers of the Gauss map.

\subsubsection{Cones}\label{Subsub:cones}

\begin{example}\label{E:cones}\rm
Let $\psi:\mathbb P^4 \dashrightarrow \mathbb P^3$ be a linear projection with center at  $p\in \mathbb P^4$. The fibers of $\psi$ define a foliation by lines $\G$ passing through $p$. Let $\f$ be the pull-back by $\psi$ of some foliation on $\mathbb P^3$, that is, $\f=\psi^*(\eta)$. It is easy to see that the leaves of $\G$ are contained in the fibers of $\GF$. Hence $\f$ is a degenerate foliation on $\mathbb P^4$. If $\eta$ is not degenerate, then $\G$ coincides with the foliation determined by the fibers of $\GF$ and $\f$ has rank three.  

Since $\psi$ is a rational map, we notice that $\f$ has a rational first integral if and only if $\eta$ has a rational first integral. The leaves of $\f$ are cones over the leaves of the foliation $\eta$ on $\mathbb P^3$.
\end{example}

\subsubsection{Joins}\label{subsub:joins}

\begin{example}\label{E:joins}\rm
Let $\G$ be a foliation by lines on $\mathbb P^4$ determined by the family of lines which intersect an irreducible surface $S$ and an irreducible curve $C$ (distinct of a line). Foliations by lines on $\mathbb P^4$ with the above property were classified in \cite{Mar,Dep3}. In general if $S$ is not a plane then $C$ is a planar curve (cf. \cite[p.396]{Mar} or \cite[Lemma 2]{Dep}). 

It follows from Lemma \ref{L:det} that
\begin{eqnarray*}
H_{\G}=E_1+2E_2\;,\;\;S=f(E_1)\;\;\text{and} \; \;C=f(E_2).
\end{eqnarray*}
We can consider the projection along the lines of $\G$ (see Definition \ref{D:projalong})
\begin{eqnarray*}
\phi_{E_1}:\mathbb P^4 \dashrightarrow S. 
\end{eqnarray*}

Given some foliation $\eta$ on $S$, we define 
\[
\f:=\phi_{E_1}^*(\eta).
\]
The leaves of $\f$ are joins $J(C,N)$, that is, an union of lines which intersect $C$ and $N$  where $N$ is a leaf of $\eta$. Thanks to Terracini's Lemma one obtain that the tangent space of $\f$ along an invariant line  by $\G$ is constant, i.e. $\G$ is tangent to the foliation determined by the fibers of $\GF$. Therefore $\f$ is a degenerate foliation on $\mathbb P^4$. 

It follows from Proposition \ref{P:rankthree} that if $\f$ does not have a rational first integral then it has rank three because  $\G$ is not tangent to a foliation by planes containing a line. Since $\phi_{E_1}$ is a rational map, $\f$ has a rational first integral if and only if $\eta$ has a rational first integral. 
\end{example}

\subsubsection{Bands}\label{subsub:bands}

\begin{example}\rm\label{E:curva}
Let $C\subset \mathbb P^4$ be a planar curve distinct from a line and $\Sigma\cong\mathbb P^2\subset \mathbb P^4$ the plane containing $C$. 

Let us define the foliation by lines. We first consider a non-constant morphism 
\[
\psi:\Sigma^*\cong\mathbb P^1 \longrightarrow C. 
\]
For each $\pi\in\Sigma^*$ let $B_{\pi}\subset \mathbb G(1,4)$ the family of lines contained in $\pi$ passing through $\psi(\pi)\in C$. Let $\G$ be the foliation by lines on $\mathbb P^4$ determined by the family
\[
B=\overline{\cup_{\pi\in\Sigma^*}B_{\pi}} \subset \mathbb G(1,4).
\]

Given $p\in C$ a smooth point, one denote by $l_p$ the tangent line to $C$ at $p$. We consider the family of planes defined in the following 
\begin{eqnarray*}
\mathcal P = \overline{\cup_{\pi}\cup_{\xi\in {(l_{\psi(\pi)})}_{\pi}^*}\xi} \subset \mathbb G(2,4),
\end{eqnarray*}
where $\pi \in {\Sigma}^*$ runs over points such that $\psi(\pi)$ is a smooth point to $C$.

Now we can define the surface
\begin{eqnarray*}
X=\left\{(\pi,\xi) \in \Sigma^*\times \mathcal P \,|\, \xi \subset \pi\right\}\subset \mathbb P^1\times\mathbb G(2,4).
\end{eqnarray*}
and consider the natural projection
\begin{eqnarray*}
\phi_X:\mathbb P^4 &\dashrightarrow& X \\
   p\in \xi &\mapsto& (\pi,\xi)\,; \;\; \xi\subset\pi.
\end{eqnarray*}

Finally, given some foliation $\eta$ on $X$ one may define 
\[
\f:=\phi_X^*(\eta). 
\]
According to the construction, the leaves of $\f$ are bands over $C$. A simple local computation shows that the tangent space of $\f$ along an invariant line  by $\G$ is constant, i.e. $\G$ is tangent to the foliation determined by the fibers of $\GF$. Then $\f$ is a degenerate foliation on $\mathbb P^4$. 

It follows from Proposition \ref{P:rankthree} that if $\f$ does not have a rational first integral then it has rank three because $\G$ is not tangent to a foliation by planes containing a line. Since $\phi_{X}$ is a rational map, $\f$ has a rational first integral if and only if $\eta$ has a rational first integral. 
\end{example}

\begin{example}\rm\label{E:riccati}
Let $\Sigma\cong\mathbb P^2\subset\mathbb P^4$ be a plane and  
\[
\psi:\Sigma^*\cong\mathbb P^1\longrightarrow \Sigma^{\lor}\cong\mathbb P^2
\]
a non-constant morphism, where $\Sigma^{\lor}$ denote the set of lines in $\Sigma$. 

We consider the surface
\begin{eqnarray*}
X=\left\{(\pi,p)\in \Sigma^*\times\Sigma \,|\, p\in \psi(\pi)\subset \Sigma\right\}\subset \mathbb P^1 \times \mathbb P^2
\end{eqnarray*} 
with natural projections
\begin{eqnarray*}
\tau_1:X &\longrightarrow& \Sigma^* \\
    (\pi,p) &\mapsto& \pi ,
\end{eqnarray*}
\begin{eqnarray*}
\tau_2:X &\longrightarrow& \Sigma \\
      (\pi,p) &\mapsto& p.
\end{eqnarray*}

Let $R$ be the divisor of ramification points of $\tau_2$.
We can write $R=D+V$, where $D$ is the divisor associated with the horizontal components (i.e. if $E$ is a component of the support of $D$, then  ${\tau_1}_{|_E}$ is dominant) and $V$ is the divisor associated with the vertical components (i.e. ${\tau_1}|_{|V|}$ is not dominant).

Now, let us define the foliation by lines. We consider a morphism 
\[
\varphi:\Sigma^*\longrightarrow\mathbb G(1,4) 
\]
satisfying
\begin{enumerate}
\item $\varphi(\pi)\subset\pi$, for all $\pi\in\Sigma^*$;
\item $\varphi(\pi)\cap\Sigma = p_{\pi}\notin \psi(\pi)$, for a generic $\pi\in\Sigma^*$.

\end{enumerate}  
where we still denote by $\varphi(\pi)$ the corresponding subset of $\mathbb P^4$. 

If $p_{\pi}$ does not belong to $\psi(\pi)$, then the family of lines 
\[
B_{\pi}=\left\{l \in \mathbb G(1,4)\;|\; l \; \text{intersect} \; \varphi(\pi) \; \text{and} \; \psi(\pi)\right\} 
\]
define a foliation by lines on $\pi$. Therefore we may consider the foliation by lines $\G$ on $\mathbb P^4$ determined by the family of lines
\begin{eqnarray*}
B=\overline{\cup_{\pi}B_{\pi}} \subset \mathbb G(1,4),
\end{eqnarray*}
where $\pi$ runs over the points of $\Sigma^*$ such that $p_{\pi}\notin \psi(\pi)$.

The foliation $\G$ defines a foliation $\eta$ on $X$ in following form. Let $\theta_{\pi}$ be the radial foliation on $\Sigma$ determined by lines passing through $p_{\pi}$ and $u_{\pi}\in {\rm{H}}^0(\Sigma,T\Sigma\otimes T_{\theta_{\pi}}^*)$ a tangent field to $\theta_{\pi}$. Since $\tau_2$ is generically a local isomorphism, for all $\pi\in\Sigma^*$ we  can define a field of directions on $\tau_2^{-1}(\psi(\pi))$
\[
v_{\pi}:={\tau_2}^*({u_{\pi}}_{|\psi(\pi)}).
\]
Hence we obtain a foliation $\eta$ on $X$ which is tangent to $v_{\pi}$ when $\pi$ runs over $\Sigma^*$. We shall say that $\G$ induces $\eta$ on $X$. 

It is not hard to prove the following proposition which gives a better knowledge about $\eta$. See \cite[Chapter 4]{Bru} for Riccati foliations.
\begin{proposition}\label{P:riccati}
Let $\G$ be the foliation by lines defined as above, then  the foliation $\eta$ induced by $\G$ on $X$ is a Riccati foliation with respect to $\tau_1$  tangent to the kernel of $d\tau_2$ along the support of $D$. Reciprocally, each Riccati foliation $\eta$ on $X$ (with respect to $\tau_1$) tangent to the kernel of $d\tau_2$ along the support of $D$, defines a foliation by lines $\G$ on $\mathbb P^4$ as above which induces $\eta$ on $X$.
\end{proposition}

Let us define the degenerate foliation on $\mathbb P^4$. Since $\Sigma$ is a fundamental component of $\G$, there exist a horizontal component $E\subset\Lambda$ such that $\Sigma=f(E)$. By using the projection along the lines of $\G$
\begin{eqnarray*}
\phi_E:\mathbb P^4 \dashrightarrow \Sigma,
\end{eqnarray*}
we obtain the rational map
\begin{eqnarray*}
\phi_X:\mathbb P^4 &\dashrightarrow& X \\
        p\in\pi          &\mapsto&   (\pi,\phi_{E}(p)). 
\end{eqnarray*}

Finally, we define the degenerate foliation on $\mathbb P^4$ as
\[
\f:=\phi_X^*(\eta). 
\]
By construction, one leaf of $\f$ is a band over a curve $N\subset \Sigma$ where $N$ is the direct image by $\tau_2$ of a leaf of $\eta$. A simple local account shows that the tangent space of $\f$ along an invariant line  by $\G$ is constant, i.e. $\G$ is tangent to the foliation determined by the fibers of $\GF$. Then $\f$ is a degenerate foliation on $\mathbb P^4$. 

It follows from Proposition \ref{P:rankthree} that if $\f$ does not have a rational first integral then it has rank three because $\G$ is not tangent to a foliation by planes containing a line. Since $\phi_{X}$ is a rational map, $\f$ has a rational first integral if and only if $\eta$ has a rational first integral. 

\begin{remark}\rm
When ${\rm{Im}}(\psi) \subset \Sigma^{\lor}$ is a line, then $|D|$ is invariant by $\eta$ and $\f$ is pull-back by $\phi_E$ of some foliation on $\Sigma$.
\end{remark}
\end{example}

\begin{example}\rm\label{E:bernoulli}
Let $\Sigma\cong\mathbb P^2 \subset \mathbb P^4$ be a plane and  
\[
\psi:\Sigma^*\cong\mathbb P^1 \longrightarrow |\mathcal O_{\Sigma}(d)|,
\]
$d > 1$ (if $d=1$ we are in Example \ref{E:riccati}),  a non-constant morphism such that  there exist a morphism
\[
\mu:\Sigma^* \longrightarrow \Sigma
\]
in which $\mu(\pi)\in\psi(\pi)$ and $\mu(\pi)$ is a point of multiplicity $(d-1)$ for $\psi(\pi)$ (where we still denote by $\psi(\pi)\subset \Sigma$ the corresponding subset of $\mathbb P^4$). 

We consider the surface 
\begin{eqnarray*}
X=\left\{(\pi,p)\in \Sigma^*\times\Sigma \,|\, p\in \psi(\pi)\right\}\subset \mathbb P^1 \times \mathbb P^2
\end{eqnarray*} 
with natural projections
\begin{eqnarray*}
\tau_1:X &\longrightarrow& \Sigma^* \\
      (\pi,p)  &\mapsto&  \pi,
\end{eqnarray*}
\begin{eqnarray*}
\tau_2:X &\longrightarrow& \Sigma \\
   (\pi,p) &\mapsto& p.
\end{eqnarray*}

Let $\theta_{\pi}$ be the radial foliation on $\Sigma$ determined by the lines passing through $\mu(\pi)$ and $u_{\pi}\in {\rm{H}}^0(\Sigma,T\Sigma\otimes T_{\theta_{\pi}}^*)$ a tangent field to $\theta_{\pi}$. Since $\tau_2$ is generically a local isomorphism, for all $\pi\in\Sigma^*$ we  can define a field of directions on $\tau_2^{-1}(\psi(\pi))$
\[
v_{\pi}:={\tau_2}^*({u_{\pi}}_{|\psi(\pi)}).
\]
Therefore we obtain a foliation $\eta$ on $X$ which is tangent to $v_{\pi}$ when $\pi$ runs over $\Sigma^*$. This case is distinct of the Example \ref{E:riccati} in the sense that the foliation $\eta$ is determined by the morphism $\psi$. 

Let us define the foliation by lines. We consider a morphism 
\[
\varphi:\Sigma^* \longrightarrow \mathbb G(1,4) 
\]
such that
\begin{enumerate}
\item $\varphi(\pi)\subset\pi$ for all $\pi\in\Sigma^*$;
\item $\varphi(\pi)\cap \Sigma=\mu(\pi)$ for all $\pi\in\Sigma^*$.

\end{enumerate}  
where we still denote by $\varphi(\pi)$ see as subset of $\mathbb P^4$. The family of lines  $B_{\pi}\subset \mathbb G(1,4)$ which intersects  $\varphi(\pi)$ and $\psi(\pi)$ defines a foliation by lines on $\pi$. Therefore one can consider the foliation by lines $\G$ on $\mathbb P^4$ determined by the family 
\begin{eqnarray*}
B=\overline{\cup_{\pi\in \Sigma^*}B_{\pi}}.
\end{eqnarray*}

Since $\Sigma$ is a fundamental component of $\G$, there exist a horizontal component $E\subset\Lambda$ such that $\Sigma=f(E)$. By using the projection along the lines of $\G$
\begin{eqnarray*}
\phi_E:\mathbb P^4 \dashrightarrow \Sigma,
\end{eqnarray*}
we obtain the rational map
\begin{eqnarray*}
\phi_X:\mathbb P^4 &\dashrightarrow& X \\
        p\in\pi          &\mapsto&   (\pi,\phi_{E}(p)). 
\end{eqnarray*}

Finally, we define the degenerate foliation on $\mathbb P^4$ as
\[
\f:=\phi_X^*(\eta). 
\]
In the construction, the generic leaf of $\f$ is a band over a curve $N\subset \Sigma$ where $N$ is the direct image by $\tau_2$ of a leaf of $\eta$. A simple local computation shows that the tangent space of $\f$ along an invariant line  by $\G$ is constant, i.e. $\G$ is tangent to the foliation determined by the fibers of $\GF$. Then $\f$ is a degenerate foliation on $\mathbb P^4$. 

It follows from Proposition \ref{P:rankthree} that if $\f$ does not have a rational first integral then it has rank three because $\G$ is not tangent to a foliation by planes containing a line. Since $\phi_{X}$ is a rational map, $\f$ has a rational first integral if and only if $\eta$ has a rational first integral. 
\end{example}

\begin{example}\rm\label{E:veronese}
Let $S=S_{1,2}$ the rational normal scroll of degree three in $\mathbb P^4$, that is, a linear projection from Veronese surface $V\subset\mathbb P^5$ with center at $p\in V$
\[
 \xymatrix{
 V  \ar@{<-}[d]_{v}  \ar@{-->}[dr]^{\tau_p}   \\
  \mathbb P^2 \ar@{-->}[r]^{\xi} & S .}
\]

We consider the linear system of dimension two, $\mathcal C:=\xi_*|{\mathcal O_{\mathbb P^2}(1)}|$ which a generic element is an irreducible conic in $S$. It is known the family of planes
\[
\mathcal P = \left\{\pi \in \mathbb G(2,4)\;|\; C \subset \pi \; \text{for some}\; C \in \mathcal C \right\},
\]
defines a foliation by planes on $\mathbb P^4$. We denote by $\pi_C$ the plane containing $C$.

We will use such foliation by planes to define a foliation by lines on $\mathbb P^4$. Let $\theta$ be a foliation on  $\check{\mathbb P}^2$ with non-degenerate Gauss map    
\begin{eqnarray*}
G_{\theta}:\check{\mathbb P}^2 &\dashrightarrow& \mathbb P^2\\
               l &\mapsto& G_{\theta}(l).      
\end{eqnarray*}
If we see $l$ as subset of $\mathbb P^2$, then $G_{\theta}(l)\in l$.

Therefore one has the following commutative diagram  
\[
\xymatrix { \mathcal C=\xi_*|{\mathcal O_{\mathbb P^2}(1)}| \ar@{<-}[d]_{\xi_*} \ar@{-->}[rr]^{\psi_{\theta}} && S  \ar@{<--}[d]_{\xi}\\ \check{\mathbb P}^2=|{\mathcal O_{\mathbb P^2}(1)}|  \ar@{-->}[rr]^{G_{\theta}} && \mathbb P^2} 
\]
where $\psi_{\theta}:=\xi\circ G_{\theta}\circ\xi_*^{-1}$. If we see $C$ as subset of $S$, then $\psi_{\theta}(C)\in C$ because $G_{\theta}(l)\in l$.

Let us denote by $\psi_{\theta}(C)_{\pi_{C}}^*\subset \mathbb G(1,4)$ the family of lines contained in $\pi_C$ passing through $\psi_{\theta}(C)$. Hence, the family  
\begin{eqnarray*}
B:=\overline{\cup_{C\in\mathcal C}\, \psi_{\theta}(C)_{\pi_{C}}^*}
\end{eqnarray*}
defines a foliation by lines $\G$ on $\mathbb P^4$ such that $\G|_{\pi_C}$ is the radial foliation on $\pi_C$ determined by the lines passing through $\psi_{\theta}(C) \in C$.

To define the degenerate foliation on $\mathbb P^4$ one consider the rational map
\begin{eqnarray*}
\varphi:\mathbb P^4 &\dashrightarrow& \check{\mathbb P}^2\\
          q\in\pi_C &\mapsto& \xi_*^{-1}(C).
\end{eqnarray*}

Therefore, we define the degenerate foliation on $\mathbb P^4$ as
\[
\f:=\varphi^*(\theta).
\]

Now, let us look at the leaves of $\f$. Since $G_{\theta}$ is generically a local isomorphism, then it sends $\theta$  to a web $\theta^{\lor}$ on $\mathbb P^2$ (the dual web of $\theta$). Hence, we obtain a web $\xi_*\theta^{\lor}$ on $S$ because $\xi$ is a birational map. We will see that the leaves of $\f$ are locally a band over a curve invariant by $\xi_*\theta^{\lor}$.

Let $p\in S$ a generic point and $C_1,...,C_k \in \mathcal C$ such that 
\[
\psi_{\theta}^{-1}(p)=\left\{C_1,...,C_k\right\}.
\]
We take $l_i\in\check{\mathbb P}^2$, $i=1,...,k$, such that $\xi_*(l_i)=C_i$. Since $G_{\theta}$ is generically a local isomorphism,  there exist one neighbourhood $U\subset S$ of $p$ and one neighbourhood $V_i \subset \check{\mathbb P}^2$ of $l_i$, $i=1,...k$, such that
\begin{enumerate}
\item[(a)] We have $\xi\circ G_{\theta}(V_i)=U$ for all $i=1,...,k$; 
\item[(b)] The map $G_{\theta}$ send $\theta_{|V_i}$ to a smooth foliation $\eta_i$ on $U$ for all $i=1,...,k$;
\item[(c)] The foliations $\eta_1,...,\eta_k$ are the local decomposition of  $(\xi_*\theta^{\lor})$ in $U$, i.e. the  tangent spaces of $\eta_1,...,\eta_k$ are tangent to $(\xi_*\theta^{\lor})$ and are in general position at all point of $U$;
\item[(d)] The foliation $\eta_i$ is tangent to $C_i$ at $p$ for all $i=1,...,k$.

\end{enumerate}

Let us consider the projection along the leaves of $\G$
\begin{eqnarray*}
\phi:\mathbb P^4 &\dashrightarrow& S \\
             p\in\pi_C &\mapsto \psi_{\theta}(C).
\end{eqnarray*}
We notice that
\begin{eqnarray*}
\phi^{-1}(p)=\cup_{i=1}^k\pi_{C_i}  
\end{eqnarray*} 
and
\begin{eqnarray*}
\phi^{-1}(U)=\cup_{i=1}^{k}W_i
\end{eqnarray*}
with $\pi_{C_i} \subset W_i$. 

It follows from the commutativity of the diagram
\[
 \xymatrix{
 \mathbb P^4\ar@/^0.5cm/@{-->}[rr]^{\phi}  \ar@{-->}[dr]_{\varphi} \ar@{-->}[r]  & \xi_*|{\mathcal O}_{\mathbb P^2}(1)|  \ar@{<-}[d]^{\xi_*} \ar@{-->}[r]_{\psi_{\theta}}  & S \ar@{<--}[d]^{\xi}
  \\
   & \check{\mathbb P}^2 \ar@{-->}[r]_{G_{\theta}} & \mathbb P^2 }
\]
that $\f|_{W_i}= \phi^*(\eta_i)$, $i=1,...,k$. We notice that since $\eta_i$ is tangent to $C_i$ at $p$, one obtains that $\pi_{C_i}$ contains the tangent line to $\eta_i$ at $p$. Therefore, by the construction, a leaf of  $\f|_{W_i}$ is a band over a leaf of $\eta_i$ because is a union of such planes $\pi_{C_i}$. Thus $\f$ is a degenerate foliation on $\mathbb P^4$.

It follows from Proposition \ref{P:rankthree} that if $\f$ does not have a rational first integral then it has rank three because  $\G$ is not tangent to a foliation by planes containing a line. Since $\varphi$ is a rational map, $\f$ has a rational first integral if and only if $\theta$ has a rational first integral. 
\end{example}

\subsection{Classification Theorem}
\begin{theorem}\label{T:teor_classif}
Let $\f$ be a degenerate foliation on $\mathbb P^4$. Then we have one of the following possibilities:
\begin{enumerate}
\item $\f$ has a rational first integral;
\item $\f$ is one of the foliations as in examples of \S \ref{Su:theexamples}.

\end{enumerate}
\end{theorem}

\begin{proof}
Let $\f$ be a degenerate foliation on $\mathbb P^4$. If $\GF$ has rank at most two, then either $\f$ has a rational first integral or $\f$ is a linear pull-back of some foliation on $\mathbb P^2$ (cf. \cite{CLn}). Since a linear pull-back of some foliation on $\mathbb P^2$ is also a linear pull-back of some foliation on $\mathbb P^3$ as in Example \ref{E:cones}, it remains to consider the case in which $\f$ has rank three.

Let $\f$ be a foliation of rank three on $\mathbb P^4$.  We can consider the foliation by lines $\G$ determined by the fibers of $\GF$. Let $M_0$ be the saturation of a germ of  leaf of  $\f$ by $\G$, as in \S \ref{Sub:focalpoints} one has the following diagram
\[
 \xymatrix{
\Lambda_{M_0}\ar@/^0.5cm/@{->}[rr]^{f_{M_0}:=f\circ i_1}   \ar[d]_{g_{M_0}} \ar[r]_{i_1}  & \Lambda   \ar[d]_{g} \ar[r]_{f}  & \mathbb P^4
  \\
  B_{M_0} \ar[r]_{i_2} & B }.
\]
Given $L\in B$ we still denote by $L$ the correspondent line see as subset of $\mathbb P^4$. Given $L$ an invariant line by $\G$  not contained in $\sing(\G)$ and seeing  the focal points as divisors in $L$, it follows from Theorem \ref{T:focais} that: 
\begin{enumerate}
\item $\Delta_{M_0}(L)$ has degree 2;
\item $\Delta_{\G}(L)$ has degree 3;
\item $\Delta_{M_0}(L)\le \Delta_{\G}(L)$. 
\end{enumerate}

\medskip

Let $E$ be an irreducible component of the support of $H_{\G}$, on one hand we can consider the codimension one foliation  $\f_E:=(f^*\f)|_{E}$ on $E$. Since the line $g^{-1}(L)$ is invariant by $f^*\f$, thus $E$ is not invariant by $f^*\f$. On the other hand we can consider the foliation $\eta_E$ defined by $f$ on $E$ (see Definition \ref{D:projalong}). We denote by $\eta_E\subset \f_E$ when $\eta_E$ is tangent to $\f_E$. The following lemma relates such foliations with the focal points.

\begin{lemma}\label{L:classifi}
Let $\f$ be a foliation of rank three on $\mathbb P^4$ and $\G$ the foliation by lines determined by the fibers of $\GF$. If $E$ is an irreducible component of the support of $H_{\G}$ such that $f(E)$ has dimension two, then  
\begin{eqnarray*}
\eta_E\subset \f_E
\end{eqnarray*} 
if and only if $|\Delta_{M_0}^*(L)|\cap E$ is not empty for a generic line invariant by $\G$.
\end{lemma}

\begin{proof}
It follows from Proposition \ref{P:depoi} that
\begin{eqnarray*}
{\rm{ker}}df(x) \subset T_xE, 
\end{eqnarray*}
for a generic point $x\in E$. Thus
\[
{\rm{ker}}df(x)={\rm{ker}}\left(df(x)_{|T_xE}\right)=T_x\eta_E.
\]

By definition $x$ belongs to $|\Delta_{M_0}^*(L)|$ if and only if $f_{M_0}$ is not an immersion at $x$. Therefore $x$ belongs to $|\Delta_{M_0}^*(L)|$ if and only if
\[
T_x\eta_E={\rm{ker}}df(x)\subset T_x(M_0\cap E)=T_x\f_E.
\] 
This proves the lemma.
\end{proof}

The Lemma \ref{L:classifi} will play a main role. Because if $g^{-1}(L)\cap E$ is a set of a single point we may consider the projection along the lines of $\G$
\[
\phi_E:\mathbb P^4 \dashrightarrow f(E),
\]
and it is easy to see that a fiber of $\phi_E$ is an invariant cone  by $\f$ if and only if
\begin{eqnarray*}
\eta_E\subset \f_E.
\end{eqnarray*} 

\medskip

If $\G$ is the radial foliation determined by the lines passing through $p\in\mathbb P^4$, then $\f$ is the linear pull-back of some foliation on $\mathbb P^3$. In fact, let $\psi:\mathbb P^4 \dashrightarrow \mathbb P^3$ be the linear projection from $p$. By hypothesis $\psi$ has fibers which is  tangent to $\f$, thus it follows from \cite[Lemma 2.2]{CLnLPT} (applied to $\psi_{|\left(\mathbb P^4 \backslash \{p\}\right)}$) that there is a foliation $\eta$ on $\mathbb P^3$ such that $\f=\psi^*(\eta)$. Therefore we are in \underline{Example \ref{E:cones}}. From now on we will suppose $\G$ does not have isolate singularities.

\medskip

Since $\Delta_{M_0}(L)$ has degree two, we will separate in the following cases:
\begin{itemize}
\item[\textbf{(I)}] - $\Delta_{M_0}(L)=p_1+p_2$, with $p_1 \neq p_2$. 
\item[\textbf{(II)}] - $\Delta_{M_0}(L)=2p$. 
\end{itemize} 

\medskip

\paragraph{\textbf{Case (I)}}
To deal with this let us consider the following:
\begin{enumerate}
\item[\textbf{(I.a)}] The support of $H_{\G}$ has just one irreducible component;
\item[\textbf{(I.b)}] The support of $H_{\G}$ has at least two irreducible distinct components.

\end{enumerate}

\medskip

\paragraph{\textbf{Case (I.a)}}
We suppose $H_{\G}=E$. The inequality 
\begin{eqnarray*}
\Delta_{M_0}(L)\le \Delta_{\G}(L)
\end{eqnarray*}
implies that $|\Delta_{M_0}^*(L)|\cap E$ is not empty for a generic line $L$ invariant by $\G$. It follows from \ref{L:det} that $f(E)$ has dimension two. Hence by Lemma  \ref{L:classifi} we get 
\[
\eta_E\subset \f_E. 
\]

Since the multiplicity at a generic point of $E$ must be constant we have that
\[
\Delta_{\G}(L)=p_1+p_2+p_3.  
\]
But this contradicts the fact that $\Delta_{M_0}(L)$ has degree two, because $\eta_E\subset \f_E$ implies that $p_i\in |\Delta_{M_0}(L)|$, $i=1,2,3$. 
Therefore this case cannot happen.

\medskip

\paragraph{\textbf{Case (I.b)}}
We first consider the following claim:
\begin{claim}
Suppose that $\Delta_{M_0}(L)=p_1+p_2$, the support of $H_{\G}$ has at least two irreducible distinct components and $\Delta_{\G}$ has pure dimension two. Then $\f$ has a rational first integral.
\end{claim} 
\begin{proof}
It follows from hypothesis
\[
\Delta_{M_0}(L)=p_1+p_2
\]
and from Lemma \ref{L:classifi} that there are $E_1,E_2$ irreducible components of the support of $H_{\G}$ (with possibility $E_1=E_2$) such that
\begin{eqnarray*}
\eta_{E_i}\subset \f_{E_i} \,,\;\; i=1,2.
\end{eqnarray*}
Thus
\[
H_{M_0}=F_1+F_2, 
\]
where $F_1$,$F_2$ are irreducible distinct components of $\Lambda_{M_0}$ and invariant by $\eta_{E_1}$, $\eta_{E_2}$ respectively. Therefore $M_0$ is a join $J(f(F_1),f(F_2))$, that is, the set of lines which intersect $f(F_1)$ and $f(F_2)$. 

Let $N_i$, $i=1,2$, be a leaf of $\eta_{E_i}$ such that  $N_i \cap F_i\ne \emptyset$. We set $i\in\{1,2\}$, since the leaves of $\eta_i$ are algebraic curves one obtain that $g^{-1}(g(N_i))$ intersect $F_j$, $i\ne j$ in a algebraic curve $C_j$ which is not invariant by  $\eta_{E_j}$. Hence $M_0$ is contained in join $J(f(C_1),f(C_2))$ which is algebraic. Therefore a generic leaf of  $\f$ is algebraic. Since $\f$ has infinite algebraic leaves, the result follows from \cite[Theorem 3.3]{Jou}. 
\end{proof}

To finish the Case (I) it remains to consider the case in which the support of $H_{\G}$ has at least two irreducible distinct components and $\Delta_{\G}$ has one fundamental component of dimension one, that is,
\begin{eqnarray*}
H_{\G}=E_1+2E_2\,,\;\; {\rm{dim}}f(E_1)=2 \;\;\text{and}\;\;{\rm{dim}}f(E_2)=1.
\end{eqnarray*}
Since $\Delta_{M_0}(L)=p_1+p_2$, it follows from Lemma \ref{L:classifi} that $\eta_{E_1}\subset \f_{E_1}$. Let us consider the projection along the lines of $\G$
\[
\phi_{E_1}: \mathbb P^4 \dashrightarrow f(E_1).
\]
The condition $\eta_{E_1}\subset \f_{E_1}$ implies that the fibers of $\phi_{E_1}$ are invariant cones by $\f$. We notice also that the fibers of $\phi_{E_1}$ are connected because these are cones over the irreducible curve $f(E_2)$. It follows from \cite[Lemma 2.2]{CLnLPT} that there is a foliation $\eta$ on $f(E_1)$ such that 
\[
\f=\phi_{E_1}^*(\eta). 
\]
Therefore $\f$ is one of the foliations as in \underline{Example \ref{E:joins}}. This concludes Case (I).

\medskip

\paragraph{\textbf{Case (II)}}
As we have noticed, a generic leaf of $\f$ is either a cone, a join or a band. By hypothesis a leaf of $\f$ cannot be a join. Suppose that a generic leaf is a cone. Then it must be a cone with vertex running over a curve $C$ (contained in the fundamental set of $\G$) because we are supposing $\G$ does not have isolate singularities.  We can consider the projection along the lines of $\G$
\[
\phi:\mathbb P^4 \dashrightarrow C,
\]
According to the hypothesis, $\f$ must be the foliation defined by the fibers of $\phi$. Therefore $\f$ has a rational first integral.  

Let us suppose that a generic leaf of $\f$ is a band, that is, a union of tangents planes to a curve. The family of such planes $\mathcal P \subset \mathbb G(2,4)$ defines a foliation by planes on $\mathbb P^4$. By Theorem \ref{T:class_folh_planos} one can separate in the following cases:

\medskip

\paragraph{\textbf{Case (1) of Theorem \ref{T:class_folh_planos}}:} If there is a line $l$ such that $\mathcal P$ is determined by the family of planes which contains  $l$, then
\[
{\rm{Im}}(\GF)\subset l^*\cong\mathbb P^2\subset \check{\mathbb P}^4. 
\]
Thus $\f$ does not have rank three. Therefore this case cannot happen.

\medskip

\paragraph{\textbf{Case (2) of Theorem \ref{T:class_folh_planos}}:} 
Suppose there exists a plane $\Sigma$ and a non-constant morphism
\[
\psi:\Sigma^*\cong\mathbb P^1 \longrightarrow \mathbb G(1,4)
\]
such that $\mathcal P$ is determined by the family of planes
\begin{eqnarray*}
\overline{\cup_{\pi\in\Sigma^*}\cup_{\xi\in(\psi(\pi))_{\pi}^*}\xi}.
\end{eqnarray*} 

Since $\mathcal P$ is tangent to $\G$, for each $\pi\in\Sigma^*$ there exist a morphism  
\begin{eqnarray*}
\varphi_{\pi}:(\psi(\pi))_{\pi}^* \longrightarrow \pi\,, \;\; \varphi_{\pi}(\xi)\in \xi,
\end{eqnarray*}
such that $\G|_{\xi}$ is determined by the family of lines which contains the point $\varphi_{\pi}(\xi)$. We have the following possibilities:

\medskip

\paragraph{\textbf{(i)}} $\psi(\pi)$ is contained in $\Sigma$ for all $\pi\in \Sigma^*$ (see Figure $2$).

\medskip

\paragraph{\textbf{(i.a)}} If ${\rm{Im}}\varphi_{\pi} = \psi(\pi)$, then
\begin{eqnarray*}
H_{\G}=3E\,,\;\;f(E)=\Sigma.
\end{eqnarray*} 
We can suppose $C^{\lor}:={\rm{Im}}(\psi)$  is a curve contained in $\Sigma^{\lor}$. Let $C \subset \Sigma$ be the dual curve of $C^{\lor}$, i.e. the family of lines  $\psi(\pi)\subset \Sigma$ are tangent to $C$. 

We set $p\in\mathbb P^4\backslash\sing(\f)$ and $M_0$ the saturation of a germ of leaf of $\f$ by $\G$ passing through  $p$. We know that $M_0$ is a band over a curve $C_{M_0}$, that is, $M_0$ is a union of planes which contain the tangent lines to $C_{M_0}$. We also know that the restriction of $\G$ to any of this planes is radial by a point $q\in C_{M_0}$. Thus
\[
C_{M_0}^{\lor}\subset C^{\lor}\subset \Sigma^{\lor},
\]
and hence $C_{M_0}\subset C$. But this contradicts the hypothesis ${\rm{Im}}\varphi_{\pi} = \psi(\pi)$. Therefore this case cannot happen.

\medskip

\paragraph{\textbf{(i.b)}} If ${\rm{Im}}\varphi_{\pi} = \{p_{\pi}\}\subset\psi(\pi)$, then $\G$ is one of the foliations as in Example \ref{E:curva}, i.e.
\begin{eqnarray*}
H_{\G}=3E\,,\;\;{\rm{dim}}f(E)=1 \;\; \text{and} \;\; f(E)\subset\Sigma.
\end{eqnarray*} 
Hence the projection 
\[
\phi_X:\mathbb P^4 \dashrightarrow X
\]
as in Example \ref{E:curva} has connected fibers tangent to  $\f$. It follows from \cite[Lemma 2.2]{CLnLPT} that there exist a foliation  $\eta$ on $X$ such that $\f=\phi_X^*(\eta)$. Therefore $\f$ is one of the foliations as in \underline{Example \ref{E:curva}}.

\medskip

\paragraph{\textbf{(i.c)}} Suppose that  
${\rm{Im}}\varphi_{\pi}$ is not contained in $\psi(\pi)$. We consider the surface
\begin{eqnarray*}
S=\overline{\cup_{\pi\in\Sigma^*}{\rm{Im}}\varphi_{\pi}} \subset \mathbb P^4.
\end{eqnarray*}
We set $p\in\mathbb P^4$ and $M_0$ the saturation of a germ of  leaf of $\f$ by $\G$ passing through  $p$. We know that $M_0$ is a band over a curve $C_{M_0}$, i.e. $M_0$ is a union of planes which contain the tangent lines to $C_{M_0}$. We also know that $\G$ restrict to any of this planes is radial by a point $q\in C_{M_0}$. Therefore $C_{M_0}\subset S$.

Since $p$ is a generic point, we may suppose that $M_0$ is not contained in $\pi$ because $\f$ is not the pencil of hyperplanes containing $\Sigma$. Thus $C_{M_0}$ is not contained in $\pi$ and the same holds for $T_qC_{M_0}$. 

But by  hypothesis under $M_0$, there exist a plane $\xi \in \mathcal P$  such that $\G$ restrict to $\xi$ is the radial foliation through $q$ and also $\xi$ contains $T_qC_{M_0}$. Since $\xi\subset \pi$ this is a contradiction, because $T_qC_{M_0}$ is not contained in $\pi$. Therefore this case cannot happen.

\medskip

\begin{figure}[h]
\centering
\includegraphics[width=.9\textwidth]{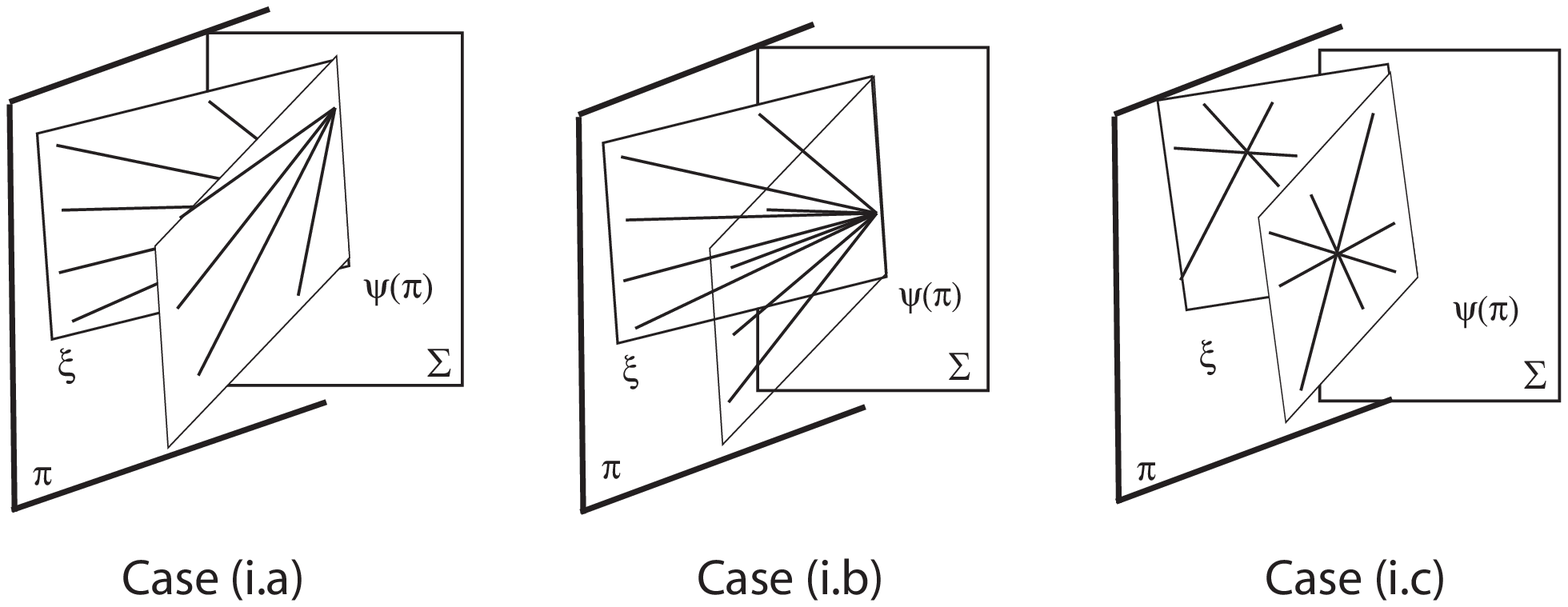}
\caption{Case (i).}
\end{figure}

\paragraph{\textbf{(ii)}} $\psi(\pi)$ is not contained in $\Sigma$ for generic  $\pi\in \Sigma^*$  (see Figure $3$). 

\medskip

\paragraph{\textbf{(ii.a)}} If ${\rm{Im}}\varphi_{\pi} = \psi(\pi)$, we consider the surface
\begin{eqnarray*}
S=\overline{\cup_{\pi \in \Sigma^*}{\rm{Im}}\varphi_{\pi}} \subset \mathbb P^4.
\end{eqnarray*}
With the same argument as in (i.c) we see that this case cannot happen.

\medskip

\paragraph{\textbf{(ii.b)}} If ${\rm{Im}}\varphi_{\pi}=\{p_{\pi}\} \subset \psi(\pi)$, we consider the curve
\begin{eqnarray*}
C=\overline{\cup_{\pi\in\Sigma^*}p_{\pi}}. 
\end{eqnarray*}
By the description of leaves of $\f$, the planes in $\mathcal P$ contains the tangent lines to $C$. But this is a contradiction because such planes are tangent to the pencil of hyperplanes  $\pi \in \Sigma^*$. Therefore this case cannot happen.

\medskip

\paragraph{\textbf{(ii.c)}} If ${\rm{Im}}\varphi_{\pi} \subset \pi\backslash(\psi(\pi)\cup \Sigma)$, we consider the surfaces
\begin{eqnarray*}
S_1=\overline{\cup_{\pi\in\Sigma^*}\psi(\pi)} \;\;\;\text{and}\;\;\;S_2=\overline{\cup_{\pi\in\Sigma^*}{\rm{Im}}\varphi_{\pi}}.
\end{eqnarray*}

We notice that 
\[
\Delta_{\G}(L)=p_1+p_2+p_3
\]
where $p_1=L\cap S_1$, $p_2=L\cap L_2$ and $p_3=L\cap \Sigma$. But this contradicts the hypothesis $\Delta_{M_0}(L)=2p$ because $\Delta_{M_0}(L)\le\Delta_{\G}(L)$. Therefore this case cannot happen.

\medskip

\paragraph{\textbf{(ii.d)}} If ${\rm{Im}}\varphi_{\pi} \subset \Sigma$, we consider the surface
\begin{eqnarray*}
S=\overline{\cup_{\pi\in\Sigma^*}\psi(\pi)}.
\end{eqnarray*} 

We suppose first that the family of curves ${\rm{Im}}\varphi_{\pi}$, $\pi\in\Sigma^*$ is constant, that is, the map $\pi \mapsto {\rm{Im}}\varphi_{\pi}$ is constant. We know that if $M_0$ is the saturation of a germ of  leaf of $\f$ by $\G$ passing through  $p$, then it is a band over a curve $C_{M_0}$, i.e. $M_0$ is a union of planes which contain the tangent lines to $C_{M_0}$. We have also that $\G$ restrict to any of this planes is radial by a point $q\in C_{M_0}$. Thus $C_{M_0}\subset {\rm{Im}}\varphi_{\pi} \subset \Sigma$. But this is a contradiction because a generic $\xi \in \mathcal P$ is not tangent to ${\rm{Im}}\varphi_{\pi}$. Therefore we may suppose that $\pi \mapsto {\rm{Im}}\varphi_{\pi}$ is not constant.

If every curve ${\rm{Im}}\varphi_{\pi}$, $\pi \in \Sigma^*$ is a line, then $\G$ is one of the foliation as in Example \ref{E:riccati}, if not,  $\G$ is one of the foliation as in Example \ref{E:bernoulli}. In both cases we have the projection 
\[
\phi_X:\mathbb P^4 \dashrightarrow X
\]
which fibers are connected and invariant by $\f$. It follows from \cite[Lemma 2.2]{CLnLPT}  that there is a foliation $\eta$ on $X$ such that $\f=\phi_X^*(\eta)$. Therefore $\f$ is either one of the foliations as in \underline{Example \ref{E:riccati}} or as in \underline{Example \ref{E:bernoulli}}.

\medskip

\begin{figure}[h]
\centering
\includegraphics[width=1.0\textwidth]{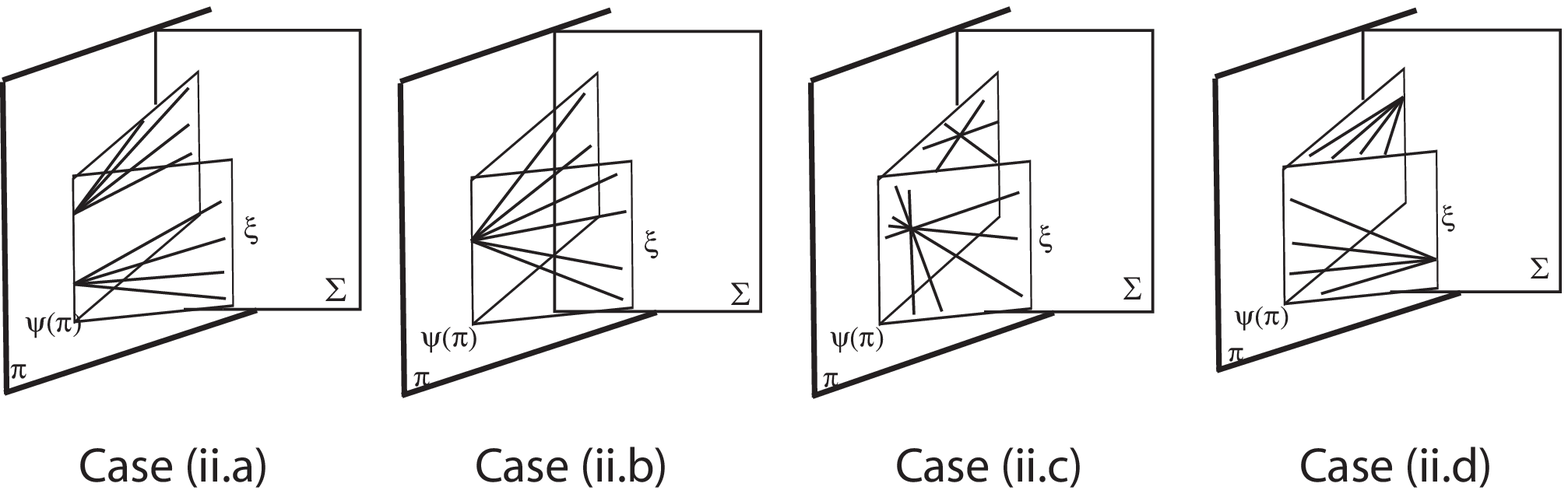}
\caption{Case (ii).}
\end{figure}

\paragraph{\textbf{Case (3) of Theorem \ref{T:class_folh_planos}}:} 
There is a cubic surface $S$ and a family $\mathcal C$ of dimension two of conics in $S$ (distinct from double lines) such that $\mathcal P$ is determined by the family of planes containing such conics. It follows from a classical Segre's Theorem (cf. \cite{S3}; \cite[Theorem $4$]{MP}) that either $S$ is a cone if such conics are reducible  or  $S=S_{1,2}$ (the rational normal scroll of degree three in $\mathbb P^4$) if a generic member of $\mathcal C$ are irreducible. 

We denote by $\pi_C$ the plane containing $C\in \mathcal C$. Since  $\pi_C$ is invariant by $\G$ for all $C\in \mathcal C$ (by definition of $\mathcal P$) one obtains that $\G$ restrict to $\pi_C$  is the radial foliation determined by the lines passing through $q_C \in \pi_C$. 

We claim that $q_C\in S$. In fact, if $q_C \notin S$ we would have 
\[
\Delta_{\G}(L)=q_C+p_1+p_2
\]
where $p_1, p_2 \in S$ because for a generic point in $S$ there is infinite invariant lines by $\G$ passing through it. But this contradicts the hypothesis $\Delta_{M_0}(L)=2p$ because $\Delta_{M_0}(L)\le\Delta_{\G}(L)$. 

We first suppose that $S$ is a cone. The family $\mathcal C$ defines a fibration by line in $S$ passing through the vertex. Let $M_0$ the saturation of a germ of  leaf of $\f$ by $\G$ through $p$. We know that $M_0$ is a band over a curve $C_{M_0}\in S$, i.e. $M_0$ is a union of planes $\xi\in\mathcal P$  which  contain the tangent lines to $C_{M_0}$. This implies that  $C_{M_0}$ is tangent to fibration defined by $\mathcal C$, hence $C_{M_0}\subset C$ for some $C\in \mathcal C$. Thus $M_0$ is a band over a line. In this case it is easy to see that $M_0$ has rank one, that is, $\f$ does not have rank three. Therefore this case cannot happen.

Let us suppose $S=S_{1,2}$. We consider the morphism 
\begin{eqnarray*}
\psi :\mathcal C &\longrightarrow& S \\
           C     &\mapsto&   \psi(C)=q_C\in C,
\end{eqnarray*}
and the projection
\begin{eqnarray*}
\tau:\mathbb P^4 &\dashrightarrow& \mathcal C\\
             p \in \pi_C &\mapsto& C.
\end{eqnarray*}
Hence, as in Example \ref{E:veronese}, we have the following commutative diagram
\[
 \xymatrix{
 \mathbb P^4 \ar@{-->}[dr]_{\varphi} \ar@{-->}[r]^{\tau}  & \mathcal C  \ar@{<-}[d]^{\xi_*} \ar@{-->}[r]^{\psi}  & S \ar@{<--}[d]^{\xi}
  \\
   & \check{\mathbb P}^2 \ar@{-->}[r] & \mathbb P^2 }
\]
Where $\varphi:=\xi_*^{-1}\circ\tau$. Since the fibers of $\varphi$ are connected and tangent to $\f$ it follows from \cite[Lemma 2.2]{CLnLPT}  that there is a foliation $\theta$ on $\check{\mathbb P}^2$ such that
$\f=\varphi^*(\theta)$. Therefore $\f$ is one of the foliations as in \underline{Example \ref{E:veronese}}.  This finish the proof of  Theorem \ref{T:teor_classif}.
\end{proof}


\begin{thebibliography}{}
\bibitem{AG}
M.A. Akivis, V.V. Goldberg.
\newblock {Differential geometry of varieties with degenerate Gauss maps,} CMS Books in Mathematics/Ouvrages de Mathématiques de la SMC, 18. Springer-Verlag, New York, (2004).

\bibitem{Bru}
M. Brunella.
\newblock { Birational geometry of foliations}.
\newblock Monografias de Matemática. [Mathematical Monographs]. Instituto de Matemática Pura e Aplicada (IMPA), Rio de Janeiro, (2000).

\bibitem{C}
D. Cerveau. 
\newblock {Feuilletages en droites, équations des eikonales et autres équations différentielles,}
\newblock {\em arXiv:math.DS/0505601v1} (2005)

\bibitem{CLn}
D. Cerveau, A. Lins-Neto.
\newblock {Irreducible components of the space of holomorphic foliations of degree two in $\mathbb P^n$, }
\newblock {\em Annals of Mathematics} 143, 577-612, (1996).

\bibitem{CLnLPT}
D. Cerveau, A. Lins-Neto, F. Loray, J.V. Pereira, F. Touzet.
\newblock {Algebraic reduction theorem for complex codimension one singular foliations,}
\newblock{\em  Comment. Math. Helv.}  81, n\textordmasculine 1, 157-169, (2006). 

\bibitem{Dep}
P. De Poi.
\newblock {Congruence of lines with one-dimensional focal locus,}
\newblock {\em Portugaliae Mathematica (N.S)},61, n\textordmasculine 3, 329-338, (2004).

\bibitem{Dep2}
P. De Poi.
\newblock {On first order congruences of lines in $\mathbb P^4$ with generically non-reduced fundamental surface,}
\newblock {\em arXiv:math.AG/0407341v3} ,(2007).

\bibitem{Dep3}
P. De Poi.
\newblock {On first order congruences of lines of $\mathbb P^4$ with a fundamental curve,}
\newblock {\em Manuscripta Math.},106, 101-116, (2001).

\bibitem{Dep4}
P. De Poi.
\newblock {On first order congruences of lines in $\mathbb P^4$ with irreducible fundamental surface,}
\newblock {\em Mathematische Nachrichten},278, n\textordmasculine 4, 363-378, (2005).

\bibitem{FPe}
T. Fassarella, J.V. Pereira.
\newblock{On the degree of polar transformations. An approach through logarithmic foliation,}
\newblock{\em Selecta Math. (N.S.),} 13, 239-252, (2007).

\bibitem{FP}
G. Fischer, J. Piontkowski.
\newblock{Ruled Varieties. An Introduction to Algebraic Differential Geometry,} 
\newblock{\em Advanced Lectures in Mathematics,} Friedr. Vieweg \& Sohn, Braunschweig, (2001).


\bibitem{GH}
P. Griffiths, J. Harris.
\newblock {Algebraic geometry and local differential geometry,}
\newblock {\em  Ann. Sci. École Norm. Sup. (4)  12}, no. 3, 355-452, (1979).

\bibitem{IL}
T.A. Ivey, J.M. Landsberg.
\newblock {Cartan for beginners: differential geometry via moving frames and exterior differential systems,} Graduate Studies in Mathematics, 61, American Mathematical Society, Providence, RI,  xiv+378, (2003). 


\bibitem{Jou}
J.P. Jouanolou.
\newblock {Équations de Pffaf algébriques,}.
\newblock {\em Lect. Notes in Math.}, 708,  (1979).

\bibitem{Kum}
E.E. Kummer.
\newblock {Über die algebraischen Strahlensysteme, insbesondere über die der ersten und zweiten Ordnung,} 
\newblock {Abh. K. Preuss. Akad. Wiss. Berlin,} 1-20 (1866), also in  
\newblock {E.E.Kummer. Collected Papers, Spring Verlag}, (1975).

\bibitem{Mar}
G. Marletta.
\newblock {Sui complessi di rette del primo ordine dello spazio a quattro dimensioni,}
\newblock {\em Rendiconti del Circolo Matematico di Palermo,} XXVIII, 353-399, (1909).

\bibitem{Mar2}
G. Marletta.
\newblock {Sui complessi di rette del primo ordine dello spazio a quattro dimensioni,}
\newblock {\em Rendiconti del Circolo Matematico di Palermo,} XXVIII, 353-399, (1909).

\bibitem{MP}
E. Mezzetti, D. Portelli.
\newblock{A tour through some classical theorems on algebraic surfaces}.
\newblock{\em An. St. Ovidius Constanta}, 5, n\textordmasculine 2, 51-78, (1997).

\bibitem{MT}
E. Mezzetti, O. Tommasi.
\newblock{On projective varieties of dimension n+k covered by k-spaces,}
\newblock{\em Illinois J.Math.}, 46 n\textordmasculine 2, 443-465, (2002).

\bibitem{PY}
J.V. Pereira, S. Yuzvinsky.
\newblock{Completely Reducible Hypersurfaces in a Pencil,}
\newblock{\em arXiv:math/0701312}, (2007).


\bibitem{Ran}
Z. Ran.
\newblock {Surfaces of order 1 in Grassmannians,}
\newblock {\em J. Reine Angew. Math.}, 368 ,119-126, (1986).

\bibitem{Rog}
E. Rogora.
\newblock {Classification of Bertini's series of varieties of dimension less than or equal to four,}
\newblock {\em Geom. Dedicata}, 64 , n\textordmasculine 2, 157-191, (1997).

\bibitem{S}
C. Segre.
\newblock {Preliminari di una teoria delle varietà luoghi di spazi,}
\newblock {\em Rend. Circ. Mat. Palermo XXX,} 87-121,(1910).

\bibitem{S2}
C. Segre.
\newblock {Su una classe di superficie degl'iperspazii legate colle equazioni lineari alle derivate parziali di 2\textordmasculine ordine,}
\newblock {\em Atti della R. Accademia delle Scienze di Torino XLII,} 559-591, (1906-1907).

\bibitem{S3}
C. Segre.
\newblock {Le superficie degli iperspazi con una doppia infinità di curve piane o spaziali,}
\newblock {\em Atti della R. Accademia delle Scienze di Torino 56,} 75-89, (1921).

\bibitem{ZILO}
F.L. Zak, A.V. Inshakov, S.M. Lvovski, A.A. Oblomkov. 
\newblock On congruences of lines of order one in $\mathbb P^3$, preprint.

\end{thebibliography}
\end{document}